\documentclass[11pt,twoside]{article}
\usepackage{amsmath, amsthm, amsfonts, amssymb}
\usepackage[all]{xy}
\usepackage{graphics}

\def\titulo#1{\noindent{\bf\LARGE{#1}} \bigskip \thispagestyle{plain}}
\def\autor#1{\noindent{\sc #1}\smallskip}
\def\direccion#1{\noindent #1\bigskip}
\def\email#1{\vspace{1cm}\noindent E-mail address: {\sf #1} \bigskip}

\theoremstyle{plain}
\newtheorem{thm}{Theorem}[subsection]
\newtheorem{proposition}[thm]{Proposition}
\newtheorem{corollary}[thm]{Corollary}

\theoremstyle{remark}
\newtheorem{remark}[thm]{Remark}
\newtheorem{remarks}[thm]{Remarks}
\newtheorem{generalities}[thm]{Generalities}

\newtheorem*{notation}{Notation}
\theoremstyle{definition}
\newtheorem{definition}[thm]{Definition}
\newtheorem{alternative definition}[thm]{Alternative Definition}
\newtheorem{example}[thm]{Example}

\def\R{\mathbb{R}}
\def\N{\n}

\def\set#1{\{#1\}}

\def\ga{\mathcal G pd\mathcal A tl}
\def\gpd{\mathcal G pd}
\def\ss{\mathcal{SS} et}
\def\cat{\mathcal C at}

\def\then{\Rightarrow}
\def\xto#1{\xrightarrow{#1}}
\def\mono{\hookrightarrow}
\def\epi{\twoheadrightarrow}

\def\a{\alpha}
\def\suba{_{\a}}
\def\b{\beta}
\def\subb{_{\b}}
\def\c{\gamma}
\def\subc{_{\c}}
\def\d{\delta}

\def\l{\lambda}
\def\set#1{\{#1\}}

\def\H{\mathcal H}
\def\U{\mathcal U}
\def\A{\mathcal A}
\def\Z{\mathbb{Z}}
\def\R{\mathbb{R}}
\def\N{\mathbb{N}}
\def\G{\mathcal G}
\def\I{\mathcal I}

\def\wt#1{\widetilde{#1}}
\def\-#1{\overline{#1}}

\def\ls{\leqslant}
\def\gs{\geqslant}
\def\<{\langle}
\def\>{\rangle}
\def\phiab{\phi\suba^{\b}}
\def\aA{\set{arrows\ of\ A}}
\def\pA{\set{points\ of\ A}}
\def\colim#1#2{\underset{#1}{colim}\ #2}

\begin{document}

\titulo{Classical Invariants for Global Actions and\\ Groupoid Atlases}

\autor{Mat\'ias L. del Hoyo, Elias Gabriel Minian}

\direccion{Departamento  de Matem\'atica\\
FCEyN, Universidad de Buenos Aires\\
Buenos Aires, Argentina.
}

\begin{abstract}

\noindent A global action is the algebraic analogue of a topological manifold. This construction was introduced in first place by A. Bak  as a combinatorial
approach to K-Theory and the concept  was later generalized by Bak, Brown, Minian and Porter to the notion of groupoid atlas. In this paper we define and
investigate homotopy invariants of global actions and groupoid atlases, such as the strong fundamental groupoid, the weak and strong nerves, 
classifying spaces and homology groups. We relate all these new invariants to classical constructions in topological spaces, simplicial complexes
and simplicial sets. This way we obtain new combinatorial formulations of classical and non classical results in terms of groupoid atlases.
\end{abstract}

\noindent{\small \it 2000 Mathematics Subject Classification.
\rm 19D99, 20L05, 18G55, 55U35.}

\noindent{\small \it Key words and phrases. \rm Global actions,
Groupoid atlases, Simplicial objects,
Homology, K-Theory.}

\section{Introduction}

The \it global action \rm construction of K-theory, introduced by A.Bak~\cite{Bak1,Bak2}, associates to a ring $A$ an algebraic object, namely 
the global action $GL(A)$ which constitutes the algebraic analogue of the standard topological construction. The underlying set of the global action
$GL(A)$ consists of the points of the general linear group of $A$ and the action consists of the virtual triangular subgroups of the general linear group
acting on the general linear group by left multiplication.

\smallskip

This new approach introduced by Bak  has the advantage that the solutions of  algebraic problems can be followed algebraically step by step.
The notion of global action gives algebraic objects such as groups, structures which allow one to develop  homotopy theory similarly to the 
classical way, by defining paths and deformations of morphisms.
In \cite{Min1, Min2}, the second named author developed an axiomatic homotopy theory for categories with natural cylinders, which 
can be applied to global actions.

\smallskip

Recently,  A. Bak, R. Brown, E.G. Minian and T. Porter generalized  ideas and constructions of global actions and 
introduced {\it groupoid atlases} \cite{BBMP}. As it
was pointed out in \cite{BBMP}, there were many advantages in formulating the concept of global action in terms of groupoids instead of group actions, so 
that it becomes part of a wider notion, namely the concept of a groupoid atlas. This was done using
the well know transition of group actions to groupoids. A groupoid
atlas can be regarded as an {\it algebraic manifold}, where the local groupoids play the role of the charts.

\smallskip

In this paper we define and investigate homotopy invariants of global actions and groupoid atlases. We study the strong fundamental groupoid, the 
weak and strong nerves, 
classifying spaces and homology theory of groupoid atlases. We also relate all these new invariants to classical constructions in topological spaces, simplicial complexes
and simplicial sets and obtain this way new combinatorial formulations of classical and non classical results in terms of groupoid atlases.

\smallskip

The rest of the paper is organized as follows.

\smallskip

In section \ref{s2} we recall the basic definitions, examples and results on global actions and groupoid atlases. Nothing is very new in this section with
the exception of a couple of new examples. One of these examples appears naturally when the global action $A(G,\H)$ (defined in \cite{BBMP}) acts on 
a $G$- set $X$. This induces a new and interesting global action. In the particular case that the general linear group acts on the quadratic forms, this 
construction is related to hermitian K-Theory.

It is important to remark that the model for the line $L$ that we introduce here is not exactly the same one that is used in \cite{Bak1} or \cite{BBMP}. 
Both models are isomorphic at the weak level, but the original notion of the line  becomes quite  rigid when one works with morphisms
 that also preserve information of the local actions. This change is essential when defining the strong fundamental group of a groupoid atlas.

We explain below the reasons of this change. We shall also prove  that both models are \it equivalent \rm groupoid atlases. Moreover,
 our model is the \it regularization \rm of the original one.
 
\smallskip

In section \ref{s3} we introduce some of  the fundamental concepts of this work, such as the notion of \it equivalence \rm between 
maps of groupoid atlases and the 
notions of 
\it irreducible \rm and \it regular \rm atlases.

\smallskip

In section \ref{s4} we study the strong fundamental group of a groupoid atlas. We use first a \it geometric \rm approach and later we prove that 
it can also be
 defined and computed with a more algebraic approach,  related to the vertex group of the colimit groupoid of the atlas. We relate the strong fundamental
groupoid with  the weak one (introduced in \cite{BBMP}) and we also show how to compute the fundamental groupoid of a topological space using 
groupoid atlases, via the Van Kampen Theorem. New formulations of  classical and non classical results on the fundamental group of an open covering are also
discussed.

\smallskip

The last section of the paper is devoted to simplicial methods. To each groupoid atlas we associate  a simplicial set, which we call the {\it (strong) nerve}
 for historical reasons and also because it generalizes in some sense the nerve of a groupoid. We also introduce a weak version of the nerve (compare with \cite{BBMP}).
We obtain this way another definition for the fundamental group of a groupoid atlas. This independently defined simplicial version is proved to coincide with 
the other two defined in the previous section.

We finish the paper with an introduction to the homology theory of groupoid atlases. We compute some easy but clarifying examples. The article ends 
with a result which relates the local homology groups with the homology of the whole atlas.

\bigskip

\section{Preliminaries: Global Actions and Groupoid Atlases}\label{s2}

Before recalling the basic definitions of {\it global actions}, 
let us begin with an easy example, which was already exhibited in \cite{BBMP}.

Let $G$ be a group and let 
$\H=\set{H\suba}_{\a\in\phi}$ be a family of subgroups of $G$, acting on $G$ by left multiplication. This will induce
a global action, denoted by $A(G,\H)$.

\smallskip

If $\H$ consists of a single group $H$, then this global action is simply the set of orbits of the action, but when $\H$ 
consists of more than one subgroup of $G$, then the different orbits of the actions interact. We will see later that 
this interaction is crucial from the homotopical point of view.
For example, consider the case  $A=A(D_3,\H)$, where $D_3$ denotes the dihedral group of order $6$. Take  
 $\phi=\set{a,b}$, $H_a=\<r\>$ with $r$ a rotation of order $3$ and $H_b=\<s\>$ with $s$ a symmetry.
The actions of each $H\suba$ divides $G$ in orbits which determine the following covering of $G$.

\

\centerline{\scalebox{1}{\includegraphics{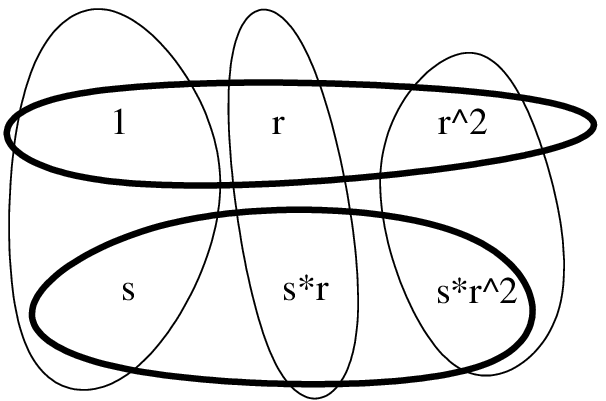}}}

\

We denote by $\U$ this covering. Its nerve $N\U$ is the simplicial complex corresponding to the following diagram
$$\xymatrix{&H_a \ar@{-}[dl] \ar@{-}[d] \ar@{-}[dr]\\
H_b & H_b.r & H_b.r^2\\
 & H_a.s  \ar@{-}[ul] \ar@{-}[u] \ar@{-}[ur]  }$$

In section 3 we shall study the strong fundamental group of global actions. As we can guess, in this case the fundamental group will be isomorphic to $\Z*\Z$. 

\subsection{Definitions and examples on global actions}

\begin{definition}
A  \it global action \rm $A$ consists of a set $X_A$ 
together with a family of
group actions 
$\set{G\suba \curvearrowright X\suba | \ \alpha\in\phi_A}$
on subsets $X\suba\subseteq X_A$. These actions are related by certain morphisms which glue them together coherently. More exactly,

$$A=(X_A,\phi_A,\set{X\suba}\suba,\set{G\suba}\suba)$$
where 
 \begin{enumerate}
    \item $X_A$ is a set,
    \item $\phi_A=(\phi_A,\leqslant)$ is an index set equipped with a reflexive relation,
    \item for each $\a\in\phi_A$, there is a subset 
	$X\suba\subset X_A$ and a \it local \rm group $G\suba$ acting on $X\suba$  and
    \item if $\a\leqslant\b$, there exists
a homomorphism $\phi\suba^{\b}:G\suba\to G\subb$
 such that $g.x=\phi\suba^{\b}(g).x$
for every  $x$ in $X\suba\cap X\subb$, $g\in G\suba$.

Moreover,  $\phi\suba^{\a}=id_{G\suba}$ and $\phi_{\c}^{\b}\circ
\phi\suba^{\c}=\phi\suba^{\b}$, whenever these compositions have sense.

\end{enumerate}
\end{definition}

Note that
 $G\suba(X\suba\cap
X\subb)\subset X\suba\cap X\subb$, i.e. the 
$\a$-orbits of  each element in the intersection are included in the intersection.

\medskip

We call $A$ a {\it single domain global action} if $X_A=X\suba$ for each $\a$. If 
$X_A=\bigcup\suba X\suba$, we say that $A$ is {\it covered}.

\bigskip

The global action $A=A(G,\H)$ of above is an example of a single domain global action.
In this case, $X_A=G$, the local actions $H\suba \curvearrowright X\suba=G$ are the subgroups of the family acting by left multiplication, and 
any two indices $\a$ and $\b$ satisfy $\a\ls\b$ if and only if $H\suba\subset H\subb$. The associated group homomorphisms
 $\phiab$ are the inclusions. 

\bigskip

Another interesting example of a single domain global action is the {\it general linear global action} $GL(n,R)$, where $R$ is
an associative ring with unit. This example was already studied in \cite{Bak1} and \cite{BBMP}, but we recall it here briefly, since it is one of the motivating examples for this theory. The homotopy groups of this global action coincide with 
the K-theory groups of the ring $R$.

Let $n\in\N$. A subset $\a$ of $\Lambda=\set{(i,j)\ |\ i\neq j,\ 1\leqslant i,j\leqslant n}$ is called  {\it closed} if  every time that it contains the pairs $(i,j)$ and $(j,k)$,
then the pair $(i,k)$ is also in $\a$.
  
Consider the poset $\phi=\set{\a\subset \Lambda\ |\ \a\ closed}$ partially ordered by inclusion.

Let $G\suba=GL(n,R)\suba$ be the subgroup of  $GL(n,R)$ generated by the matrices 
 $$\set{E_{ij}(r)|\ r\in R, (i,j)\in\a},$$
where $E_{ij}(r)$ is the matrix containing $1$ in the diagonal, $r$ in the position $(ij)$ and $0$ elsewhere.

It is not difficult to verify that a matrix  $A$ belongs to $GL(n,R)\suba$ if and only if  $A_{ij}=1$ for $i=j$ and $A_{ij}=0$ if $(i,j)\notin\a$.
 
For $\a\subset\b$, we denote $\phi\suba^{\b}$ the inclusion $GL(n,R)\suba\to GL(n,R)\subb$.
Now let $X\suba=GL(n,R)$. The subgroup $GL(n,R)\suba$ acts on  $GL(n,R)$ by left multiplication.

\bigskip

Note that the general linear global action is a particular example of the actions $A(G,\H)$ introduced before. 
If $A(G,\H)$ is considered as an atlas with discrete index set (see \cite{BBMP}), this is not longer true since in the general linear global action situation we are making some identifications between local actions.

\bigskip

\begin{example}
Let $A=A(G,\H)$ and let $G\curvearrowright X$ be a $G$-set. The global action $A\curvearrowright X$ is an extension of the action of $G$ on $X$ and
it is defined as follows. Take $X_{A\curvearrowright X}=X\suba=X$ for each $\a\in\phi_{A\curvearrowright X}=\phi_A$, 
with the action induced by the action of $G$.

As a particular case, consider  the general lineal group $GL(n,\R)$ acting over the quadratic forms in $\R^n$ by base change, i.e.
 $$(C\cdot Q)(x)=Q(C.x),\ C\in GL_n(\R),\ Q \text{ a quadratic form},\ x\in\R^n.$$
This action restricts well to the subset $X$ of positive defined quadratic forms.
If $A$ is the matrix of $Q$, then $C^tAC$ is the matrix of $C\cdot Q$.

\end{example}

\subsection{Path components and the weak fundamental group}

Consider again our first example $A=A(D_3,\H)$, where $D_3$ denotes the dihedral group of order $6$, and 
take the element
$1\in D_3$. When different elements of the local groups $H_a$ and $H_b$ act consecutively on it, we obtain a path like
the following

$$1\xrightarrow{s.}s\xrightarrow{r.}rs\xrightarrow{s.}srs\xrightarrow{r.}rsrs=1.$$

In general, the various actions of the local groups $G\suba$ interact on the intersections of the local sets $X\suba$ and
this induces a global dynamics in $A$. The elements of $X_A$ move along $X_A$ through the actions of the different local 
groups. This idea suggests definitions for paths and loops in $A$ and therefore notions for path connectedness and simply
connectedness.

\bigskip

Let $A=(X_A,\phi_A,\set{X\suba}\suba,\set{G\suba}\suba)$ be a global action and let $\alpha\in \phi_A$. An {\it $\alpha$-frame}
is a finite subset $\set{x_0,\ldots,x_n}\subset X_A$ such that for each $i$ there exists $g_i\in G\suba$ with 
$g_i\cdot x_{i-1}=x_i$. % We say that $\set{x_0,\ldots,x_n}$ is a {\it local frame} if the index $\a$ was not specified. 
A (weak) {\it path} is a finite sequence $x_0,\ldots,x_n$ such that for each $i$ the set $\set{x_{i-1},x_i}$ is a {\it local frame}, i.e. an $\alpha$-frame for some index $\a$.

Given  two elements $x,y$ of $X_A$, we say that they are in the same {\it connected component} of $A$ if there exists a path joining both points. As usual, we denoted by $\pi_0(A)$ the set of (path) components of $A$.

\bigskip

Let us compute $\pi_0(A)$ in the examples of above.

Consider first the case $A=GL(n,R)$. If $x$ and $y$ are in the same component, then there is a finite sequence of matrices $E_i\in G_i=GL(n,R)_{\alpha_i}$ such that
$$E_nE_{n-1}\dots E_1 x=y$$
which implies that they determine the same class in the quotient $GL(n,R)/ E(n,R)$. Here $E(n,R)$ denotes the subgroup
of elementary matrices.
Since every elementary matrix can be factored as a finite products of $E_i\in GL(n,R)_{\alpha_i}$, we obtain that
$$\pi_0(GL(n,R))=GL(n,R)/E(n,R)=K_1(n,R).$$

For more details, see \cite{Bak1, Bak2, Bass, BBMP, Mil1, Swan}.

In the case $A=A(G,\H)$, a similar argument shows that $x,y\in A$ are in the same component if and only if there exists
a sequence $h_{\alpha_i}\in H_{\alpha_i}$ such that $h_{\a_n}\dots h_{\a_1}x=y$. If we denote
$\<\H\>=\<H_i\ |\ i\in \phi\>$ the subgroup of $G$ generated by all $H_i$, then we obtain
$$\pi_0(A(G,\H))=G/\<\H\>.$$
 
\bigskip
 
Recall that a {\it weak morphism} $f:A\to B$ is a set theoretic function $f:X_A\to X_B$ which preserves local frames. A path is a particular 
example of a weak morphism between global actions. More exactly, let $L$ be the global action with underlying set $X_L=\Z$, with the index 
set $\phi_L\subset \mathcal{P}(\Z)$ the family of subset of $\Z$ of the form $\set{n}$ and  $\set{n,n+1}$, and whose local actions are the free 
and transitive actions of the trivial group and $G_2$ respectively. A path in $A$ is simply a morphism $L\to A$ that {\it stabilizes} in both directions (i.e. $\exists\ N$ such that $f(n)=f(n+1)$ for $|n|>N$).

If a weak morphism $f:L\to A$ does not stabilize, we call it a {\it weak curve}. A path with the same initial and final point is a {\it weak loop}.

\begin{remark}
As we pointed out in the introduction, the model for the line $L$ that we introduce here is not the same one used in \cite{Bak1} or \cite{BBMP}. 
Both models are isomorphic at the weak level, but the original notion  becomes very  rigid when one works with morphisms
 that  preserve the information of the local actions.

\end{remark}
\bigskip

Recall  that the product of  global actions is defined as follows. Given two global actions $A$ and $B$, the {\it product} $A\times B$ is the global action with underlying set $X_A\times X_B$ and index set $\phi_A\times\phi_B$, equipped with the product relation. The local action of $A\times B$ indexed by $(\a,\b)$ is the product action between $G\suba\curvearrowright (X_A)\suba$ and $G\subb\curvearrowright (X_B)\subb$.

The product $A\times B$, defined as above, satisfies the universal property of the categorical product in the category of global actions.

\

A homotopy between paths $\omega$ and $\omega'$ is defined as a weak morphism
$$H:L\times L\to A$$
for which there exist integers $N_0, N_1$ such that $H(-,N_0)=\omega$ and $H(-,N_1)=\omega'$, and that stabilizes in an appropriate sense (see \cite{Bak1}). In particular, the local frames of the product $L\times L$ are the subsets $S$ satisfying
$$S\subset \set{(n,m),(n,m+1),(n+1,m),(n+1,m+1)}$$
for some $n,m\in\Z$.

\begin{definition}\label{weak pi1}
Let $\omega$ and $\omega'$ be loops based on $x\in X_A$. A {\it homotopy} between $\omega$ and $\omega'$ is a function $H:\Z\times\Z\to X_A$ such that
\begin{itemize}
\item for all $m,n$, the sets $\set{H(n,m),H(n,m+1),H(n+1,m),H(n+1,m+1)}$ are local frames of $A$, and
\item there exist $N_0<N_1,M_0<M_1\in\Z$ such that $H(-,M_0)=\omega$, $H(-,M_1)=\omega'$ and
$$H(n,m)=\begin{cases} H(N_0,m) & if\ n<N_0 \\H(N_1,m) & if\ n>N_1\\H(n,M_0) & if\  m<M_0 \\H(n,M_1) & if\ m>M_1\end{cases}$$
\end{itemize}
\end{definition}

We can thus define the weak fundamental group of a global action $X_A$ with base point $a\in X_A$ as follows.

\begin{definition}
The {\it weak fundamental group} $\pi_1^w(A,x)$ is the set of homotopy classes of loops at $x$. The multiplication is defined via concatenation of paths.
\end{definition}

Let us compute $\pi_1^w(A,x)$ in the case $A=A(G,\H)$. Without loss of generality, we may suppose $x=e$ the neutral element of $G$. We denote with $\underset{\cap}{\coprod}H\suba$ the amalgamated product of $H\suba$, i.e. the colimit in the category of groups over the diagram $\set{H\suba\cap H\subb,H\suba}$.

\begin{proposition}
The weak fundamental group of $A(G,\H)$ can be computed as the kernel of the canonical map $\underset{\cap}{\coprod}H\suba\to G$.
\end{proposition}

\begin{proof}
The group $\pi_1^w(A,e)$ is isomorphic to the fundamental group of the nerve of the cover of $G$ by $H$-coclasses, with $H\in\H$ (see below). Let $\wt G=\underset{\cap}{\coprod}H\suba$ 
and write $j(H)$ for the image of $H\to\wt G$. Let $N\wt G$ be the nerve of the cover of $\wt G$ by $j(H)$-coclasses. 
The canonical group homomorphism $\wt G\to G$ induces the universal covering $N\wt G\to NG$ between the nerves (see \cite{AH},\cite{BBMP}). 
Given $j(H)$ a vertex of $N\wt G$, a Deck transformation $\varphi$ is determined by $\varphi(j(H))$, which could be any element of $\set{g.j(H)|g\in Ker(\wt G\to G)}$, the fiber over $H$. Multiplication by $g$ gives a simplicial map on $N\wt G$ for $g\in\wt G$, hence it gives a Deck transformation for $g\in Ker(\wt G\to G)$. We conclude that the group of Deck transformations is exactly $Ker(\wt G\to G)$.
\end{proof}

{\bf Note.} Here, by a {\it covering} of simplicial complexes, we mean a simplicial map with the unique lifting property of simplices. Observe that a simplicial map $K\to L$ is a covering if and only if $|K|\to |L|$ is a covering of topological spaces.

\

At the weak level, a global action is the same as a set equipped with a cover.
Let $A$ be a global action, $X_A$ the underlying set and $\U_A$ the covering determined by the local orbits. Let $V(\U_A)$ be the {\it Vietoris complex} of $\U_A$, whose simplices are the finite subsets of $X_A$ that are included in some element of $\U_A$. Let $N(\U_A)$ be the {\it nerve} of this covering, whose simplices are the finite subsets of $\U_A$ with non trivial intersection. Since local frames in $A$ are just simplices in $V\U_A$, we have
$$\pi_0(A)\cong\pi_0(V\U_A)\cong\pi_0(N\U_A)$$
and
$$\pi_1^w(A)\cong\pi_1(V\U_A)\cong\pi_1(N\U_A).$$
Dowker's theorem, proved in \cite{BBMP}, relates both simplicial complexes $N\U_A$ and $V\U_A$.

In general, one has the following result.
\begin{proposition}\label{ser debil no sirve}
The functor $A\mapsto VA$ is an equivalence between the category of (covered) global actions with weak morphism and the category of simplicial complex.
\end{proposition}

An inverse can be obtained by giving to a simplicial complex $K$ a global action $$\set{S(s)\curvearrowright s\ |\ s \ simplex\ of\ K},$$
with the simplices of $K$ as indices and whose underlying set is the set of vertices of $K$. Here $S(X)$ means the group of bijections of the set $X$.

\

There exists a stronger notion of morphisms of global actions, which originally were called regular in $\cite{Bak1}$.

A {\it regular morphism} $f:A\to B$ is a triple $(X_f,\phi_f,G_f)$ satisfying
\begin{enumerate}
\item $\phi_f:\phi_A\to \phi_B$ is a relation preserving map,
\item $G_f(\a):G\suba\to G_{\phi_f(\a)}$ is a group morphism such that if $\a\ls\b$ the diagram
$$\xymatrix{G\suba \ar[r] \ar[d] & G_{\phi_f(\a)} \ar[d]\\
            G\subb \ar[r]        &   G_{\phi_f(\b)}}$$
commutes,
\item $X_f:X_A\to X_B$ is a set function such that $X_f(X\suba)\subset X_{\phi_f(\a)}$ and

\item for each $\a\in\phi_A$, $(G_f,X_f):G\suba\curvearrowright X\suba\to G_{\phi_f(\a)}\curvearrowright X_{\phi_f(\a)}$ is a morphism of actions, i.e. $G_f(\a)(g)\cdot X_f(x)=X_f(g\cdot x)\ \forall\ g\in G\suba,\ x\in X\suba$.
\end{enumerate}

In order to generalize the constructions of above, we must be very careful about what a strong path or a strong homotopy is.
Regular morphisms, that a priori allow us to work with curves, paths and loops in a strong sense, are very restrictive.
 To give an idea, let us consider a global action such that all its local groups are finite of odd order. 
Because a regular morphism contains group morphisms as part of its data, there are no regular morphisms $L\to A$ which are not constant. 
The regular maps $L\to A$ do not measure, in general, the dynamics of $A$.

\bigskip

{\bf Important Note.}
In this paper, the word \bf regular \rm will mean a different concept (see \ref{regu}). 
Since regular morphisms of global actions (in the sense of \cite{Bak1}) are not used in this article,  there should be no confusion.

\subsection{Groupoid atlases}

Giving a group $G$ acting on a set $X$, one can associate to the group action $G\curvearrowright X$ a groupoid $\G=G\ltimes X$. The objects of $\G$ are the elements of $X$ and for any $x,y\in X$, the arrows from $x$ to $y$ are the pairs $(g,x)$ with $g\in G$ such that $g\cdot x=y$. Composition is defined in the obvious way.

Applying this construction to each local orbit of a global action $A$, we obtain a groupoid atlas. 
The concept of a groupoid atlas was introduced in $\cite{BBMP}$. 

\begin{definition}
A {\it groupoid atlas} $A=(X_A,\phi_A,\G_A)$ consists of a set $X_A$, an index set $\phi_A$ equipped with a reflexive relation, and for each $\a\in\phi$, a groupoid $\G\suba$ such that
\begin{enumerate}
\item[-] $X\suba=Obj\ \G\suba\subset X_A$,
\item[-] if $\a\ls \b$ in $\phi_A$, $X\suba\cap X\subb$ is union of components of $\G\suba$, and there is a functor $\phiab:\G\suba|_{X\suba\cap X\subb}\to \G\subb|_{X\suba\cap X\subb}$ which restricts to the identity in objects.
\end{enumerate}
The structural functor $\phiab$ is the identity when $\a=\b$, and it commutes with the compositions. The groupoids $\G\suba$ are called the {\it local groupoids} of $A$.
\end{definition}

\begin{notation}
Let $\G$ be a groupoid and let $X\subset Obj\ \G$. We denote by $\G|_X$  the full subgroupoid of $\G$ on $X$.
\end{notation}

\begin{example}
As we mentioned above, every global action induces a groupoid atlas. In particular, we endow $A(G,\H)$ and $GL(n,R)$ with a structure of global action. Although not every groupoid atlas comes from a global action. For some examples, see \cite{Por1}.
\end{example}

\begin{example}
Any groupoid $\G$ can be viewed as a groupoid atlas with trivial index set. 
 This groupoid atlas is denoted by $a(\G)$. This construction induces  a fully faithful functor into any of the categories of groupoid atlases that 
we  discuss later.
\end{example}

\begin{example}\label{sc example}
Let $K$ be a simplicial complex. We define a groupoid atlas $a(K)$ as follows: $X_{aK}=V_K$, the set of vertices of $K$;
$\phi_{aK}=S_K$, the set of simplices of $K$ ordered by inclusion; for each simplex $s$, the local groupoid $\G_s$ is the 
simply connected groupoid (tree) with object set $s$; the structural functors $\phi_s^t$ are the inclusions.

As a particular case, the {\it n-sphere} is defined to be the groupoid atlas $a(\partial\Delta[n])$, where $\partial\Delta[n]$ is the simplicial complex with vertices $\set{0,...,n}$ and simplices all the nonempty proper subsets of $ \set{0,...,n}$.
\end{example}

This functorial construction $K\mapsto aK$ associates to every simplicial complex $K$ a groupoid atlas $aK$ satisfying some extra properties that we will discuss later. For example, $aK$ is irreducible, regular and infimum. Also, this construction induces a fully faithful functor into the category $[\ga]$ defined in the next section. This way one might view groupoid atlases as generalized simplicial complexes, for which there are many others local models than the (homotopy trivial) simplices.

\bigskip

Another interesting examples of groupoid atlases arise from the fundamental groupoids of topological spaces. 

\begin{example}\label{topological example}
Let $X$ be a topological space and $\U$ an open cover of $X$. The groupoid atlas $A(X,\U)$ is defined as follows. The underlying set is $X$ and 
the index set $\phi$ is the poset $(\U,\subset)$. For each $U\in\U$, the local groupoid $\G_U$ is the fundamental groupoid $\pi_1(U)$ and the morphisms $\phi_U^V$ are induced by the inclusions $U\hookrightarrow V$.
\end{example}

In section \ref{s4} we relate both the weak and the strong version of the fundamental group of $A(X,\U)$ with $\pi_1(X)$.

\

{\it Weak morphisms} between groupoid atlases are, as one might suppose, functions between the underlying sets which preserve local frames. 
In this context, a {\it local frame} is a finite subset of some connected component of a local groupoid. This category, with groupoid atlases as 
objects and weak morphism as arrows, is canonically equivalent to those of global actions and simplicial complexes. 
We call it the {\it weak category} of groupoid atlases.

\

The notion of (strong) morphism of groupoid atlases is not as rigid as the notion of regular morphism of global actions.

\begin{definition}
A {\it morphism} $f:A\to B$ between groupoid atlases is a triple $(X_f,\phi_f,\G_f)$ satisfying
\begin{enumerate}
\item[-] $X_f:X_A\to X_B$ is a set-theoretic function,
\item[-] $\phi_f:\phi_A\to \phi_B$ is a function which preserves the relation $\leq$,
\item[-] $\G_f:\G_A\to\G_B$ is a (generalized) natural transformation of groupoid diagrams over the function $\phi_f$, which restricts to $X_f$ in the objects.

In other words, for each $\a$ a functor $\G_f(\a):\G\suba\to \G_{\phi_f(\a)}$ is given in such a way that $Obj\ \G_f(\a)=X_f|_{X\suba}$ and, if $\a\ls\b$, the diagram
$$\xymatrix{\G\suba \ar[r] \ar[d] & \G_{\phi_f(\a)} \ar[d]\\
            \G\subb \ar[r]        &   \G_{\phi_f(\b)}}$$
commutes.
\end{enumerate}
\end{definition}

%The composite of $f:A\to B$ and $g:B\to C$ is the morphism $g\circ f:A\to C$ given by $X_{g\circ f}=X_g\circ X_f$, $\phi_{g\circ f}=\phi_g\circ\phi_f$ and for each $\a$ the functor $\G_{g\circ f}(\a)=\G_g(\phi_f(\a))\circ\G_f(\a)$.
%\bigskip
We denote by $\ga$  the category of groupoid atlases with (strong) morphisms.

\bigskip

The atlas $A$ is {\it covered} if every element of $X_A$ is an object of some local groupoid. All the
groupoid atlases that we consider are assumed to be covered.

Note that, if $A$ is covered and  $f:A\to B$ is any morphism, then the function $X_f$ is locally determined by the values of the functors 
$\set{\G_f(\a)}\suba$ in objects. However, $X_f$ must be part of the data.

To illustrate this, consider the following example. Suppose $A$ is the groupoid atlas with 
 $X_A=\ast$, $\phi_A=\set{1,2}$ with the discrete order and $X_1=X_2=\ast$ and let $B$ be any groupoid atlas.
Suppose also that we are given a map $\phi_f:\phi_A\to \phi_B$ and a generalized natural transformation $\G_f:\G_A\to \G_B$. Note that
the map $\phi_f:\phi_A\to \phi_B$ picks two indices $\phi_{f(1)}$ and $\phi_{f(2)}$ of $\phi_B$, and  $\G_f:\G_A\to \G_B$
 determines two objects $x$ and $y$  of $\G_{\phi_{f(1)}}$ and $\G_{\phi_{f(2)}}$, respectively. If we take $x\neq y$, 
then $(\phi_f,\G_f)$  does not come from a morphism of atlases.

\

Under certain condition on the atlas $A$, a relation preserving map $\phi_f:\phi_A\to \phi_B$ and a 
natural family of functors $\G\suba\to \G_{\phi_f(\a)}$ do determine a map $A\to B$. We introduce the concept of a \it good \rm atlas,  
which solves this problem and plays a fundamental role in the next sections. 
This concept  is weaker than (but intimately related with) the notion of  {\it infimum}, introduced in $\cite{Bak1}$ and $\cite{BBMP}$ (cf. \ref{infimum}).

\begin{definition}\label{def of good}
Let $A$ be a groupoid atlas. We say that $A$ is {\it good} if for every $x\in X_A$ the set
$$\phi_x=\set{\a\in\phi_A\ |\ x\in X\suba}$$
has an initial element, i.e. there exists $\a_x\in\phi_x$ such that $\a_x\ls\b$ for all $\b\in\phi_x$.
\end{definition}

\begin{remark}
Suppose $A$ is a good atlas. Given a relation preserving map $\phi_f:\phi_A\to\phi_B$ and $\G_f:\G_A\to \G_B$ a natural family 
of functors over $\phi_f$, the various local functions $\set{Obj\ \G_f(\a)}_{\a\in\phi_A}$ agree in the intersections because $Obj\ \G_f(\a)(x)$ 
must be $Obj\ \G_f(\a_x)(x)$ for all $\a$, and therefore there is a well defined function $X_f:X_A\to X_B$. In this case, a morphism $f:A\to B$ can be
 regarded as a pair $(\phi_f,\G_f)$.
\end{remark}

\section{Equivalences of Maps and Irreducible Atlases}\label{s3}

In order to define the right notion of strong fundamental group we need to change first a little bit our notion of morphism. Many maps $A\to B$ seem to carry the same information, and they just differ in the indices. For example, if $A$ is such that $\phi_A=\set{0<1}$ with $\G_0$ a subgroupoid of $\G_1$ and $\phi_0^1$ is the inclusion, this atlas does not contain more information than the one coming from $\G_1$.
Moreover, given $a\in X_A$, it is necessary to identify the various identities of this element viewed as an object of the different local groupoids.
Also, we should understand $\phiab(g)$ as an \it extension \rm of the \it movement $g$.\rm

This induces the notion of equivalence between morphisms.

\subsection{Morphisms modulo equivalences}

\begin{definition}
Given $f,f':A\to B$, we say that $f$ is a {\it corestriction} of $f'$ if $\phi_f(\a)\ls \phi_{f'}(\a)$ for all $\a\in\phi_A$ and all the diagrams
$$\xymatrix{
 & \G\suba \ar[dl]_{\G_f(\a)} \ar[dr]^{\G_{f'}(\a)} \\
\G_{\phi_f(\a)} \ar[rr]_{\phi_{\phi_f(\a)}^{\phi_{f'}(\a)}} & & \G_{\phi_{f'}(\a)}}$$
are commutative. We write $f\lhd f'$.

The morphisms $f,f':A\to B$ are {\it equivalent}  if they are in the same class of the equivalence relation generated by the corestrictions. We write $f\sim f'$.
\end{definition}

\begin{remark}
It is easy to see that $f\sim f'$ implies $X_f=X_{f'}$.
\end{remark}

Note that the relation $\lhd$ is transitive. Moreover, if $f\lhd f'$, then $g\circ f\lhd g\circ f'$ and $f\circ h\lhd f'\circ h$ for any $g,h$.
It follows that, if $f,f':A\to B$ and $g,g':B\to C$ are such that $f\sim f'$ and $g\sim g'$, then $g\circ f\sim g'\circ f'$. Therefore, we 
can define a new category $[\ga]$, with groupoid atlases as objects and the classes of (strong) morphisms as maps. We denote the class of $f$ by $[f]$.

\bigskip

The notion of equivalence between morphisms loosen the dependence of $A$ from the index set $\phi_A$. Let us go back to the example at the beginning of 
this section. The information of the atlas $A$ seems to lie just in $\G_1$. If we call $B=a\G_1$ the atlas 
with this only local groupoid, the canonical arrows $f:A\to B$ and $g:B\to A$ satisfy $gf\sim id_A$ and $fg = id_B$, so they are isomorphic in $[\ga]$.

\begin{definition}
We say that $f:A\to B$ is an {\it equivalence} if there exists $g:B\to A$ such that $gf\sim id_A$ and $fg\sim id_B$, i.e. $[f]$ is an 
isomorphism in $[\ga]$. Given atlases $A$ and $B$, we say that they are {\it equivalent} if there is an equivalence between them.
\end{definition}

\

\subsection{Irreducible atlases}

Recall that the set $\phi_A$ of a groupoid atlas $(X_A,\phi_A,\G_A)$  is  equipped with a reflexive relation. %, that can be thought antisymmetric up to equivalence
If $\a,\b\in\phi_A$ are such that $\a\ls\b$, by definition $X\suba\cap X\subb$ can be written as union of connected components of the 
groupoid $\G\suba$. If in addition $\G\suba$ is connected, it follows that $X\suba\subset X\suba\cap X\subb$ and the structural
 functor $\phiab:\G\suba|_{X\suba\cap X\subb}\to \G\subb|_{X\suba\cap X\subb}$ can be regarded simply as a functor between $\G\suba$ and $\G\subb$.
 So, when each local groupoid is connected, a groupoid atlas is nothing but a particular diagram of groupoids.

\begin{definition}
We will say that a groupoid atlas $A$ is {\it irreducible} if $\G\suba$ is a connected groupoid for every index $\a$ in $\phi_A$.
\end{definition}

Given any atlas $A$, there is a natural way to associate to it an irreducible atlas $i(A)$, 
considering each component of any local groupoid as an individual local groupoid. The underlying set of $i(A)$ is the same of $A$, so $X_{i(A)}=X_A$. 
The index set of $i(A)$ must repeat each index of $A$ as many times as many components of the corresponding groupoid. Explicitly,
$$\phi_{i(A)}=\set{\ (\a,X)\ |\ \a\in\phi_A,\ X=Obj\ \G'\subset X_A,\ \G'\ component\ of\ \G\suba}$$ 
and the relation is induced from $\phi_A$: 
$$(\a,X)\ls(\b,Y) \iff \a\ls\b\ and\ X\subset Y.$$
Finally, the local groupoid $\G_{(\a,X)}$ is the component of $\G\suba$ with object set $X$, and if $(\a,X)\ls(\b,Y)$, 
the morphism $\phi_{(\a,X)}^{(\b,Y)}$ equals $\phiab$. Like in any irreducible atlas, the structural maps of $i(A)$ are defined over the whole local groupoid.\medskip

This construction is functorial, since every map $f=(X_f,\phi_f,\G_f):A\to B$  induces a new map $if:iA\to iB$ given by $X_{if}=X_f$, $\phi_{if}(\a,X)=(\phi_f(\a),Y)$ where $Y$ is the component of $\G_{\phi_f(\a)}$ that contains $\G_f(\a)(X)$, and $\G_{if}(\a,X):\G_{(\a,X)}\to\G_{\phi_{if}(\a,X)}$ is the restriction of the functor $\G_f(\a):\G\suba\to \G_{\phi_f(\a)}$.

\

The canonical map $\varphi_A:iA\to A$ that forgets the second coordinate of the indices ($\phi_{\varphi_A}(\a,X)=\a$) and such that
 $\G_{\varphi_A}(\a,X):\G\suba|_{X}\to \G\suba$ is the inclusion of each component, is natural in $A$.
It satisfies also the following universal property, whose proof is straightforward.

\begin{proposition}
Let $g:B\to A$ a map of groupoid atlases, and suppose that $B$ is irreducible. Then, there exists a unique map $h:B\to iA$ such that $g=\varphi_A\circ h$.
\end{proposition}

It is clear that $iA=A$ when the atlas $A$ is already irreducible. The functor $i$ is a right adjoint 
for the inclusion of the full subcategory of irreducible atlases into the category of atlases.

\

The map $\phi_A:iA\to A$ is a weak isomorphism, but it is not an isomorphism or an equivalence in general. However, the groupoid atlases $A$ and $iA$ are intimately related: they share all the invariants that we will study, 
such as the nerve, the homology groups and even the fundamental group.

\

\subsection{A new approach}

Given an irreducible atlas $A$ and indices $\a,\b$ and $\c$ in $\phi_A$ such that $\a\ls\b\ls\c$, then $X\suba\subset X\subb\subset X\subc$ and 
 the composition $\phi\subb^{\c}\circ\phiab$ is defined and it  agrees with $\phi\suba^{\c}$ if $\a\ls\c$ in $\phi$. Therefore, in 
an irreducible atlas we can assume that the relation $\ls$ is transitive without loss of generality.

\begin{proposition}\label{up to equivalence}
Given  an irreducible atlas $A$, the relation $\ls$ defined in its set of indices $\phi_A$ is a partial order ``up to equivalence'', i.e. there exists an atlas $B$ whose set of indices is partially ordered by $\ls$, and an equivalence $A\to B$.
\end{proposition}

\begin{proof}
Recall that the relation $\ls$ defined in $\phi_A$ is always reflexive, and since $A$ is irreducible, it can also be taken transitive. It only 
remains to make $\ls$ antisymmetric.

The atlas $B$ will be obtained from $A$ by deleting some local groupoids. We say that two indices
$\a,\b\in\phi_A$ are {\it paired} if $\a\ls\b$ and $\b\ls \a$. Clearly, been paired is an equivalence relation. We denote by $[\a]$ the paired class of $\a$.
If $\a$ and $\b$ are paired indices,
 then $X\suba=X\subb$ and the functor $\phiab$ must be an isomorphism. Let $c$ be a selector function that assigns to each paired class an element of itself. We define $B$ as follows,
\begin{enumerate}
\item[-]$X_B=X_A$,
\item[-]$\phi_B=\set{c[\a]\ |\ \a\in\phi_A}$, with the relation induced by the inclusion $\phi_B\mono\phi_A$,
\item[-]$\G_B$ is the restriction of $\G_A$ to $\phi_B$.
\end{enumerate}

The map $f:B\to A$ is the canonical inclusion. The inverse, $g:B\to A$ is given by
\begin{enumerate}
\item[-]$X_g=id:X_A\to X_B$,
\item[-]$\phi_g:\phi_A\to\phi_B$ is the map $\a\mapsto c[\a]$. If $\a\ls\b$, because $\a$ and $c[\a]$ are paired, $\b$ and $c[\b]$ are paired and $\ls$ is transitive, we have $c[\a]\ls c[\b]$,
\item[-]the functor $\G_g(\a):\G\suba\to \G_{c[\a]}$ equals $\phi\suba^{c[\a]}$ the structural functor of $A$. The square
$$\xymatrix{
\G\suba \ar[r]^{\phi\suba^{c[\a]}} \ar[d]_{\phiab}& \G_{c[\a]} \ar[d]^{\phi_{c[\a]}^{c[\b]}} \\
\G\subb \ar[r]_{\phi\subb^{c[\b]}} & \G_{c[\b]}
}$$
is commutative because of the naturality of $\G_A$.

\end{enumerate}

The composition $g\circ f$ is the identity of $B$, and $f\circ g$ is a corestriction of $id_A$. It follows that  $A$ and $B$ are equivalent.
%Observe that the construction of $B$ is not natural: it requires to specify the indices $c[\a]$.
\end{proof}

%Anyway, redefining $\phi\suba^{\c}=\phi\subb^{\c}\circ\phiab$ if necessary, the relation $\ls$ can be closed by transitivity keeping the information of $A$. In this way, we may associate to $A$ an atlas with its set of indices partially ordered. Next, we will investigate in which sense they are related.

The next definition, which is motivated by the last proposition, gives another approach to irreducible groupoid atlases. The proof of the 
equivalence between both definitions is omitted.

\begin{alternative definition}\label{alternative definition}
Given  a partially ordered set $\phi$, an {\it irreducible groupoid atlas} $A$ with index set $\phi$ is a diagram in the 
category of groupoids $\G_A:\phi\to\gpd$ such that for each $\a$ the groupoid $\G_A(\a)=\G\suba$ is connected and if $\a\ls\b$ then the induced functor
$\phiab$, is an inclusion on the objects.
\end{alternative definition}

With this definition, the underlying set $X_A$ is the union $\bigcup Obj\ \G\suba$.

A morphism  $f:A\to B$ is a pair $(\phi_f,\G_f)$, with $\phi_f:\phi_A\to\phi_B$ a map of posets and $\G_f$ a natural transformation between $\G_A$ and $\G_B\circ\phi_f$ that induces a function $X_f$.

$$\xymatrix{\phi_A \ar[rr]^{\phi_f} \ar[dr]_{\G_A}& & \phi_B \ar[dl]^{\G_B} \\ &\gpd}$$

In this context, the corestrictions arise naturally.
Given $f,f':A\to B$, the relation $f\lhd f'$ is equivalent to the existence of a natural transformation $\eta:\phi_f\then\phi_f':\phi_A\to\phi_B$ such that $\G_A\xto{\G_f'}\G_B\circ\phi_{f'}$ equals the composition $\G_A\xto{\G_f} \G_B\circ\phi_f \xto{\G_B\circ\eta} \G_B\circ\phi_{f'}$,

$$\xymatrix{\phi_A \ar@/^.6pc/[rr]^{\phi_{f'}} \ar@/_.6pc/[rr]_{\phi_f} \ar[dr]_{\G_A}& \Uparrow & \phi_B \ar[dl]^{\G_B} \\ &\gpd}$$

\

\subsection{Regular atlases}

We may think a regular atlas as a {\it well pointed} atlas. Our model for the line is regular, and this is very convenient because of the existence of sections for the projection $A\times L\to A$ (they are necessary for the homotopies).

\begin{definition}\label{regu}
An atlas $A$ is called {\it regular} if it is good and $\G_{\a_x}=\ast$ for every $x\in X_A$, where $\a_x$ is the minimum of $\phi_x$ (cf. \ref{def of good}).
\end{definition}

\begin{remark}
Given $x\in X_A$, the index $\a_x$ in the definition of regular atlas is not unique, but two of such indices 
must be mutually related by $\ls$, and therefore the uniqueness can be assumed up to equivalence (proposition \ref{up to equivalence}) 
when $A$ is irreducible.
\end{remark}
\

There is a regularization functor $A\to r(A)$ defined as follows.
$$X_{r(A)}=X_A$$ and $$\phi_{r(A)}=\phi_A\sqcup \set{\a_x}_{x\in X_A},\ \text{ with } \a_x\ls \a \text{ if and only if } x\in X\suba.$$
The local groupoids $\G_{\a_x}$ are singletons and the $\G\suba\ 's$ are as before.

\begin{example}
The line $L$ is a regular groupoid atlas. As we pointed out above, it is the regularization of the one used in \cite{Bak1} and \cite{BBMP}.
\end{example}

Sometimes one need to replace a groupoid atlas $A$ by its regularization $rA$, in order to work with paths, homotopies and morphisms in a right way.
The following proposition asserts that, when $A$ is good, this replacement does not change the equivalence class of the atlas. Then, every good atlas can be supposed to be regular up to equivalence.

\begin{proposition}\label{a sim ra}
Let $A$ be a good atlas. The canonical inclusion $A\to rA$ is an equivalence.
\end{proposition}

\begin{proof}
Let $f$ be the inclusion $A\to rA$, we define $g:rA\to A$ by sending $\phi_g(\a_x)$ to the initial element of $\phi_x$. It is easy to see that $g$ is a morphism, $gf=id_A$ and $id_{rA}\sim fg$ (in fact, $id_{rA}\lhd fg$).
\end{proof}

We will say that a morphism $f:A\to B$ between regular atlases is {\it regular} if $\phi_f(\a_x)=\a_{X_f(x)}$ for all $\a$. 
The following three results will be used to characterize the (strong) fundamental group of groupoid atlases.

\begin{proposition}
Let $f:A\to B$ be a morphism between regular atlases. Then there exists $f':A\to B$ regular such that $f\sim f'$.
\end{proposition}

\begin{proof}
Let us define $f'$ and prove that it is equivalent to $f$. Take $X_{f'}$ equals $X_f$. For each $\a$ there is at least one $x\in X_A$ such that $\a=\a_x$ is the initial element of $\phi_x$. Define $\phi_{f'}$ as the function that 
sends the initial element of $\phi_x$ to the initial element of $\phi_{X_f(x)}$, and the other indices $\a$ to $\phi_f(\a)$. The family of functors $\G_{f'}(\a):\G\suba\to\G_{f'(\a)}$ is defined as follows: if there is an $x$ such that $\a=\a_x$, $\G_{f'}(\a)$ is uniquely determinated by $X_f$; otherwise, $\G_{f'}(\a)$ is just $\G_f(\a)$. Naturality follows from the naturality of $\G_f$, and since $f'\lhd f$ the result follows.
\end{proof}

\begin{proposition}\label{r adj inc}
Let $f:A\to B$. If $B$ is good, then there exists $g:r(A)\to B$ that extends $f$. Moreover, if $g'$ is another extension, then $g\sim g'$.
\end{proposition}

\begin{proof}
It is sufficient to extend $\phi_f$ to the indices $\set{\a_x\ |\ x\in X_A}$ of $\phi_{rA}$ that are not in $\phi_A$, and this can be done sending each $\a_x$ to the minimum of $\set{\a \in\phi_B\ |\ X_f(x)\in X\suba}$. 
The morphism $g$ obtained this way is a corestriction of any other extension.
\end{proof}

When $A$ is good and $f=id_A:A\to A$, the map $g:rA\to A$ is exactly the one constructed in \ref{a sim ra}. Proposition \ref{r adj inc} proves that the regularization functor is right adjoint to the inclusion of the full subcategory of $[\ga]$ formed by the good atlases.

\begin{proposition}\label{regular proposition}
Let $f,g:A\to B$ with $A$ regular and $B$ good. If $\phi_f(\a)=\phi_g(\a)$ and $\G_f(\a)=\G_g(\a)$ for all $\a\in\phi_A$ except perhaps for the indices $\set{\a_x}$, then $f\sim g$.
\end{proposition}

\begin{proof}
Follows from the previous result.
\end{proof}

\

\section{The (Strong) Fundamental Group}\label{s4}

Along this section, we work with good groupoid atlases.

\subsection{Global points and global arrows}

From the notion of equivalence, we will obtain a definition for points and arrows of an atlas $A$ in a global sense, i.e. independent of the local groupoids.

We denote with $*$ the singleton groupoid (a single object and a single arrow), and with $\I$ the 2-point simply connected groupoid. Recall that $a:\gpd\to\ga$ is the functor that maps each groupoid $\G$ into the atlas whose unique local groupoid is $\G$.

\

A morphism $a(\ast)\to A$ is a pair $(\a,x)$, where $\a\in \phi_A$ and $x\in X\suba$ are the images of the unique index and the 
unique object respectively. A same element $x\in X_A$ gives rise to many morphisms $a(\ast)\to A$, one for each $\a$ such that $x\in X\suba$. 
On the other hand, if $\a\ls\b$, the two corresponding morphisms are equivalent. 

\begin{definition}
A {\it global point} in $A$ is the equivalence class of an arrow $a(\ast)\to A$.
$$\pA=Hom_{[\ga]}(a(*),A)$$
\end{definition}

\begin{remark}
The canonical function $\set{points\ of\ A}\to X_A$ that maps the class of $p:a(\ast)\to A$ in $X_p(\ast)$ is bijective when the atlas $A$ is good. In fact, it is surjective because $A$ is covered and injective because each $\phi_x$ has initial element.
\end{remark}

As we pointed out before, we think $\phi\suba^{\b}(g)\in\G\subb$ as an extension of the {\it movement} $g\in \G\suba$.
Let us consider the case of $A(G,\H)$, with the family $\H$ closed under (finite) intersections. The arrows of a local groupoid $H\ltimes G$ come from the 
action of the subgroup $H$ in $G$ by left multiplication. Recall that $(h,g)$ is the arrow with source $g$ and target $hg$. These arrows are identified 
with the corresponding elements of $G$, without considering which $H$ does actually contain them: if $h$ belongs to $H$ and $H'$, since $\H$ is closed under intersections we might think $(h,g)\in H\ltimes G$ and $(h,g)\in H'\ltimes G$ as extensions of $(h,g)\in (H\cap H')\ltimes G$. Thus, the arrows of $A$ are
$$\set{(h,g)|\ g\in G,\ h\in H \ for\ some\ H\in\H}$$
For a general atlas $A$, we propose the following definition.

\begin{definition}\label{def de global arrow}
A {\it global arrow} $[g]$ of $A$ is the class of some local arrow $g\in \G\suba$ by the relation generated by $g\sim \phi\suba^{\b}(g)$. 
Equivalently, a global arrow is the equivalence class of a morphism $a(\I)\to A$.
$$\aA=Hom_{[\ga]}(a(\I),A)$$
\end{definition}

Since $g\sim g'$ implies $X_g=X_{g'}$, there are well defined {\it source} and {\it target} of a global arrow $[g]$. They are $X_g(0)$ and $X_g(1)$, respectively.

\

Let $A$ be a good atlas, and $x\in X_A$. If $x\in X\suba\cap X\subb$, then $id_x\in \G\suba \sim id_x\in \G\subb$ because we can reach one from the other by functors $\phi\suba^{\b}$. We denote this global arrow as $id_x$. Observe that $[g]=id_x$ does not imply that $g$ is the identity of $x$ in some local groupoid, because the structural functors $\phiab$ may not be faithful.

\

Finally, recall that the {\it standard $n$-simplex} $\Delta[n]$ in the category of groupoids is the simply connected groupoid with object set $\set{0,...,n}$.
In particular, we have $\ast=\Delta[0]$ and $\I=\Delta[1]$. These classical definitions give us models for the simplices in $\ga$, composing with $a$. As a generalization of what we did with points and arrows, they will be used in section \ref{s5} to define the nerve of a groupoid atlas.

\

\subsection{Curves, paths and loops}

In this subsection we will interpret a curve $L\to A$ as a sequence of global arrows. This leads to a very nice result: the fundamental group of a groupoid atlas equals the fundamental group of the colimit groupoid over the set of indices.

\

In the category $\ga$, we consider the inclusions $i_n:a(\I)\to L$ given by $\phi_{i_n}(\ast)=\set{n,n+1}$ and $id=\G_{i_n}(\ast):\I\to\G_{\set{n,n+1}}.$
A morphism $\l:L\to A$ gives rise to a sequence $(\l_n)_{n\in\Z}$ defining $\l_n:a(\I)\to A$, $\l_n=\l\circ i_n$.

\

Given $\l:L\to A$, in the triple $(X_{\l},\phi_{\l},\G_{\l})$ the function $X_{\l}$ is determinated by the others, and any natural family of functors $(\phi_{\l},\G_{\l})$ induces a map $\l:L\to A$ since $L$ is regular.
 Moreover, the local groupoids of $L$ of the form $\G_{\set{n}}$ are singletons so each functor $\G_{\l}(\set{n})$ must be trivial, and each functor $\G_{\l}(\set{n,n+1})$ is given by an arrow in $\G_{\phi(\set{n,n+1})}$(the image of the arrow $\set{n\to n+1}$). 
 Thus, a morphism $\l$ determines and is determined by the following data.
\begin{itemize}
	\item $(\a_n)_n\subset \phi_A$ where $\a_n=\phi_{\l}(\set{n})$, and
	\item $(\l_n)_n$, where $\l_n=\l\circ i_n:\I\to A$ satisfies $s(\l_{n+1})=t(\l_n)$ and $\a_n,\a_{n+1}\ls \phi_{\l_n}(*)$.
\end{itemize}

\begin{remark}
A map $L\to A$ gives rise to a curve with a {\em framing} $\b$ in the vocabulary of \cite{BBMP}. Here, $\b$ is the function $n\to \a_n$. Conversely, one can get a map $L\to A$ from a weak curve $w$ with a framing $\b$, specifying local arrows $g_n$ which make $w$ a curve.
\end{remark}

%We will say that the sequence $(\a_n)_n$ is the {\em frame} of the map $\l$ (compare with \cite{BBMP}).

\begin{definition}
A {\it curve} $l$ in $A$ is an equivalence class  of a morphism $\l:L\to A$.
\end{definition}

If $(\l_n)_n$ is a sequence of local arrows such that $s(\l_{n+1})=t(\l_n)$, we can construct a morphism $\l:L\to A$ taking for each $n$ an 
index $\a_n\ls \phi_{\l_{n-1}}(*),\phi_{\l_n}(*)$, which exists when $A$ is good. Here $\phi_{\l_n}(*)$ is the index of the groupoid that contains the local arrow $\l_n$.
Since $L$ is good, by proposition \ref{regular proposition}, the class of this $\l$ does not depend on the choice of the $\a_n$.
Hence, we have established the following (partial) correspondence:
\begin{align*}
(\l_n:a(\I)\to A)_n/ s\l_{n+1}=t\l_n &\then			[\l]\text{ curve}\\
([\l_n])_n\subset \aA				&\Leftarrow	\l:L\to A
\end{align*}

If $\l\sim\l'$ then $[\l_n]=[\l]\circ[i_n]=[\l']\circ[i_n]=[\l'_n]$ for all $n$. 
We investigate now under which hypothesis the converse is also true, namely: $\l_n\sim \l'_n\ \forall n\ \then\ [\l]=[\l']$.

\begin{proposition}
Let $\l,\l':L\to A$ be such that $\l_n=\l'_n\ \forall n\neq n_0$ and $\l_{n_0}\lhd\l'_{n_0}$. Then, $\l\sim\l'$.
\end{proposition}

\begin{proof}
For each $n$, let $\a_n\in\phi_A$ be less or equal than $\phi_{\l}(\set{n})$ and $\phi_{\l'}(\set{n})$. Let $h:L\to A$ be given for the sequences $(\a_n)_n$ and $(\l_n)_n$.
It is easy to check that $h\lhd \l$ and $h\lhd \l'$.
\end{proof}

\begin{corollary}
Let $\l,\l':L\to A$. Suppose that there exists a finite subset $J\subset \Z$ with $\l_n=\l'_n$ for $n\notin J$ and $\l_n\sim \l'_n$ for $n\in J$, then $\l\sim\l'$.
\end{corollary}

\begin{proof}
It suffices to construct a finite sequence $\l=\l^0,\l^1,\dots,\l^k=\l'$ such that $\l^i$ and $\l^{i+1}$ are in the conditions of the proposition.
\end{proof}

We say that $\l:L\to A$ {\it stabilizes} to the left (resp. to the right) if there are $x\in X_A$ and $\a\in\phi_A$ such that $\l_n=id_x\in\G\suba$ for small (resp. big) enough values of $n$.

\begin{proposition}\label{equiv entre curvas}
If $\l,\l':L\to A$ stabilize in both directions and $\l_n\sim\l_{n'}$ for every $n$, then $\l\sim\l'$.
\end{proposition}
\begin{proof}
Let $(N_0,N_1)$ be a stabilization pair for $\l$ and $\l'$.
We have $\l_n=id_{x_0}\in\G_{\a_0}$ for $n\ls N_0$ and some $\a_0$ in $\phi_A$, and
$\l_n=id_{x_1}\in\G_{\a_1}$ for $n\gs N_1$ and some $\a_1$.
The map $\l$ gives an infinit sequence $(\phi_{\l}(\set{n}))_n\in\phi_A$ that not necessary stablize in any direction. However, we may suppose that 
$\phi_{\l}(\set{n})=\a_0$ for $n\ls N_0$ and
$\phi_{\l}(\set{n})=\a_1$ for $n\gs N_1$ 
without change the equivalence class of $\l$ (actually, $\l$ is a corestriction of a map satisfying this).

Similarly, let $\a'_0$ and $\a'_1$ be the indices where $\l'$ stabilizes.

Choose $\c_0,\c_1$ in $\phi_A$ such that $\c_i\ls \a_i, \a'_i$, $i=0,1$. Define $(\d_n)_n\subset\phi_A$ as 
$$\d_n=\begin{cases}\c_0 & \text{ if } n\ls N_0\\
										\phi_{\l}(\set{n}) & \text{ if }N_0<n<N_1\\
										\c_1 & \text{ if } n\gs N_1\end{cases}.$$

Define $g_n:\I\to A$ by
$$g_n=\begin{cases}id_{x_0}\in\G_{\c_0} & \text{ si } n\ls N_0\\
									 \l_n & \text{ si }N_0<n<N_1\\
									 id_{x_1}\in\G_{\c_1} & \text{ si } n\gs N_1\end{cases}.$$

Denote $\wt{\l}:L\to A$ the morphism induced by the sequences $(\d_n)_n$ and $(g_n)_n$. Clearly, $\wt{\l}\lhd \l$. Analogously, 
define $\wt{\l'}:L\to A$ such that $\wt{\l'}\lhd \l'$. By the propostion of above and since $\wt{\l}_n=\wt{\l'}_n$ for $n<N_0$ or $n>N_1$ and 
 $\wt{\l}_n=\l_n\sim\l'_n=\wt{\l'_n}$ for $N_0\ls n\ls N_1$, it follows that  $\wt{\l}\sim\wt{\l'}$. Therefore, $\l\sim\l$.
\end{proof}

If $l$ is a curve, we say that $l$ {\it stabilizes} to the left (resp. to the right) if there exists $\l:L\to A$ such that $[\l]=l$ and $\l$ stabilizes to the left (resp. to the right).

\begin{definition}
A {\it path} in $A$ is a curve which stabilizes in both directions.
\end{definition}

Note that if $l$ is a path in $A$, then there exists $\l:L\to A$, $N_0,N_1\in\Z$, $x_0,x_1\in X_A$ and $\a_0,\a_1\in\phi_A$ 
such that $l=[\l]$, $\l_n=id_{x_0}\in \G_{\a_0}$ for $n< N_0$ and $\l_n=id_{x_1}\in \G_{\a_1}$ for $n> N_1$. 
%Furthermore, $N_0$ and $N_1$ can be taken such that $[\l_{N_0}]\neq id_{x_0}$ and $[\l_{N_1}]\neq id_{x_1}$. Then there is a  well defined notion of length of a path.

\

Consider again  the correspondence between curves and sequences of global arrows. A path $l=[\l]$ is associated to the sequence $([\l_n])_n\subset \aA$ which stabilizes in identities at both sides.

Conversely, given a sequence $(g_n)_n\subset \aA$ such that $s(g_{n+1})=t(g_n)$, $g_n=id_{x_0}$ for $n<N_0$ and $g_n=id_{x_1}$ for $n>N_1$, 
one can take local arrows $(\l_n)_n$  such that $[\l_n]=g_n$ and construct the curve associated to this sequence. 
Note that this curve does not have to be a path. This can be solved choosing $\a_0,\a_1$ such that $x_0\in X_{\a_0}$, $x_1\in X_{\a_1}$ and 
taking $\l_n=id_{x_0}\in\G_{\a_0}$ and $\l_n=id_{x_1}\in\G_{\a_1}$ for $n$ small or big enough. This new curve  is
 a path, and is uniquely determinated by proposition \ref{equiv entre curvas}. Thus we have obtained the following result.

\begin{proposition}
The constructions of above are mutually inverse. They define a bijection between the set of paths of an atlas $A$ and the sequences $(g_n)_n\subset\aA$ that stablizes on identities and satisfies $s(g_{n+1})=t(g_n)$ for all $n$.
\end{proposition}

This correspondence preserves sources and targets, and concatenations of paths up to homotopy. This leads to a simplicial computation
 of the fundamental group. In fact, we will prove in the next section 
that the fundamental group of an atlas $A$ is the fundamental group of its {\it nerve} $NA$.

\

\subsection{The fundamental group}

In this subsection we introduce the (strong) fundamental group of a good groupoid atlas $A$ with a fixed base point. Later we will 
generalize this construction to any groupoid atlas. The definitions and results introduced here admit
 formulations in terms of the fundamental groupoid.

\

Recall that $L$ is regular. Its points can be represented by the inclusions
$$j_n:\ast\to L,\ \phi_{j_n}(\ast)=\set{n},\ id=\G_{j_n}(\ast):\ast\to\G_{\set{n}}.$$

The projections $p_1,p_2:L\times L\to L$ admit sections $id_L\times j_n, j_n\times id_L$ respectively. 
We will simply write $j_n$ instead of $id_L\times j_n$, the inclusion of $L$ in the $n$-th row.

\begin{definition}
A {\it homotopy} between two paths $l$ and $l'$ is a map 
	$$H:L\times L\to A$$
such that there exist integers $N_0,N_1,M_0$ and $M_1$ satisfying
\begin{itemize}
\item $H\circ j_{N_0}=\l$ with $[\l]=l$,
\item $H\circ j_{N_1}=\l'$ with $[\l']=l'$,
\item there are $\a_0,\a_1\in\phi_A$ such that for all local arrow $f$ of $L\times L$, $\G_H(f)=id_{x_0}\in\G_{\a_0}$ if the first coordinate of $s(f)$ is less than $M_0$ and $\G_H(f)=id_{x_1}\in\G_{\a_1}$ if the first coordinate of $t(f)$ is greater than $M_1$.
\end{itemize}
\end{definition}

It is clear that  homotopy of paths is  an equivalence relation.

\begin{definition}\label{def de pi1}
Let $A$ be a good atlas and $x\in X_A$. The fundamental group $\pi_1(A,x)$ is the set of homotopy classes of loops at $x$, with the operation induced by
 the concatenation of paths.
\end{definition}

\begin{proposition}\label{pi1 y sim}
If $f:A\xto{\sim}B$ is an equivalence, then $\pi_1(A,x_0)\cong\pi_1(B,X_f(x_0))$.
\end{proposition}
\begin{proof}
It is a consequence of the definition of paths and homotopies. The fundamental group is actually a well defined invariant in the category $[\ga]$. 
For an alternative proof, see section \ref{s5}.
\end{proof}

{\bf Important remark:} We extend the definition of the fundamental group to non 
necessarily good atlases, taking the fundamental group of their regularizations. By the last proposition, this new definition agree with 
the original one on good atlases.

\

A homotopy $H:l\cong l'$ can be decomposed into more elementary homotopies, considering the various  
ways to reach the upper-right point from the lower-left one in a lattice (say $[0,N]\times[0,N]\subset X_{L\times L}$), that gives rise to a sequence of paths $...,\l_i,\l_{i+1},...$ 
where two consecutive paths just differ in a single square. Note that all these paths become loops when $H$ is applied.

These elementary homotopies relate a path $([g_{M_0}],\dots,[g_{M_1}])$, viewed as a sequence of global arrows, with another $([h_{M_0}],\dots,[h_{M_1}])$, 
where $h_{i+1}\circ h_i=g_{i+1}\circ g_i\in\G\suba$ for some $i$ and some $\a$, and $h_j=g_j\in\G_{\a_j}$ for $j\neq i,i+1$. Inserting the 
loop $([g_{M_0}],\dots,[g_i\circ g_{i-1}],\dots,[g_{M_1}])$ between them, we can observe that the homotopies in $A$ define the 
same equivalence relation that the one used to represent the morphisms in a colimit on the category of categories (cf. \cite{GZ},\cite{McL}).

\smallskip

This proves the following result.

\begin{thm} \label{fundamental group}
Let $A$ be an irreducible groupoid atlas and let $x\in X_A$. Consider 
$\G(A)=\colim{}{\G\suba}$ the colimit over the diagram $\phi_A$. Then the
fundamental group $\pi_1(A,x)$ equals the fundamental group of the groupoid $\pi_1(\G(A),x)$.
\end{thm}

\begin{remark}
The hypotesis of irreducible on $A$ is necessary to make $\set{\G\suba}\suba$ a groupoid diagram. Otherwise, the functors $\phiab$ are partially defined and the colimit does not have sense.
\end{remark}

\

Last theorem shows another way to define the group $\pi_1(A,x_0)$ when $A$ is irreducible. When $A$ is not irreducible, the group $\pi_1(iA,x_0)$ defined in this way is isomorphic to the fundamental group defined as before. This approach will be study in section \ref{s5}.

\

\begin{example}\label{topological example 2}
We return to example \ref{topological example}. By definition, the local groupoid $\G_U=\pi_1(U)$ is connected when the open subset $U$ is path connected. So, the atlas $A(X,\U)$ is irreducible if and only if $U$ is path connected for all $U\in\U$.
If in addition the family $\U$ is closed under finite intersections, then
	$$\pi_1(A(X,\U))=\colim{U\in\U}{\pi_1(U)}=\pi_1(X)$$
as a consequence of the last theorem and (the groupoid formulation) of  Van Kampen's theorem (cf. \cite{Brown}, \cite{May}).
\end{example}

\

\subsection{Some examples}

The fundamental group $\pi_1(A,x)$ appears as a natural invariant for groupoid atlases in many ways: from paths and homotopies of 
paths (\ref{def de pi1}), by simplicial calculation (\ref{pi1 de na}) and even as a purely algebraic invariant (\ref{fundamental group}, \ref{colimit groupoid}). 
In some special cases coincides with the weaker one $\pi_1^w(A,x)$. The group $\pi_1^w(A,x)$ describes some of the algebra of the global 
nontrivial loops, while $\pi_1(A,x)$  detects also the local ones. It is clear that these groups are not equal in general.

\

The following examples are the simplest cases where the differences between these two groups are noticed. In a proposition of below we will 
prove that these are essentially the unique kinds of differences they could have.

\medskip

Given  a groupoid atlas $A$ and $x\in X_A$, recall that the group $\pi_1^w(A,x)$ is the fundamental group of the simplicial complex $V\U_A$. We denote this complex by $VA$.

\begin{example}
Let $\G$ be a groupoid and let $x$ be an object of $\G$. The groups $\pi_1(a(\G),x)$ and $\pi_1^w(a(\G),x)$ depends only of the 
component of $\G$ that contains $x$, so $\G$ can be supposed connected. The Vietoris complex $V(a(\G))$ is the simplex spanned by $O=Obj\ \G$, so 
it is contractible. Hence, $\pi_1^w(a(\G),x)=\pi_1(V(a(\G)),x)=0$. On the other hand, $\pi_1(a(\G),x)=Hom_{\G}(x,x)$, which is not trivial in general. 
\end{example}

\begin{example}
Let $A$ be a groupoid atlas with discrete set of indices $\phi_A=\set{\a,\b}$ and such that $X_A=X\suba=X\subb$ and the 
local groupoids $\G\suba$ and $\G\subb$ are simply connected. Every finite subset of $X_A$ is a local frame, so as in the last example, 
$VA$ is contractible and $\pi_1^w(A,x)=0$ for any $x$. Let $y\xto{i}z$ be the unique arrow in $\G_i$ with source $y$ and target $z$.
The (strong) loop whose associated sequence of global arrows is 
$$([x\xto{\a}y],[y\xto{\b}x])$$
is a non trivial element of $\pi_1(A,x)$ for any $y\in X_A$, $y\neq x$ (in fact, $\pi_1(A,x)$ is the free group with one generator $([x\xto{\a}y],[y\xto{\b}x])$ for each $y$). Therefore, the groups $\pi_1$ and $\pi_1^w$ are distinct if $A$ has other points than $x$.
\end{example}

\

 There is a canonical group morphism
$$p:\pi_1(A,x)\to \pi_1^w(A,x)$$
which sends a loop $l$ to the weak loop $X_{\l}$, where $\l$ is a representative of $l$ that stabilizes in both directions. The map $p$ can also be 
defined using the simplicial information: mapping a path of global arrows $(g_1,...,g_N)$ to the path $(\set{s(g_1),t(g_1)},...,\set{s(g_N),t(g_N)})$ 
of edges of $VA$. This map is defined up to homotopy. Note that $p$ is the group morphism induced by the simplicial map $p$ of \ref{p en nervios}.

The last examples show that $p$ is not an isomorphism in general, but it is so in many interesting cases.

\

Consider for example the groupoid atlas $A=A(G,\H)$, with $G$ a group and $\H$ a family of subgroups closed 
under finite intersections. Let $x$ be an element of $G$. In order to compute the group $\pi_1(A,x)$, note that $\aA=\set{(h,g)|\ g\in G,\ h\in H \ for\ some\ H\in\H}$, as one can easily check from definition \ref{def de global arrow} and the paragraph above it. Then, we have the isomorphism
$$\pi_1(A,x)\xto{\sim}\pi_1^w(A,x)$$
since (the groupoid form of) $p$ maps bijectively global arrows into edges and elementary homotopy triangles into three-elements local frames. This result is 
not true if the family $\H$ is not closed under intersections.

\

The general linear global action is a particular case of an $A(G,\H)$ where  $\H$ is closed under intersections. To see that, recall that for $A=GL(n,R)$ 
the index set $\phi_A$ consists of  the closed subsets of $\set{(i,j)\ /\ i\neq j,\ 1\leqslant i,j\leqslant n}$ partially ordered by inclusion. If $\a$ and $\b$ are in $\phi_A$, it is clear that $\a\cap\b$ is also closed, and hence it is in $\phi_A$. The subgroup $GL(n,R)_{\a\cap\b}$ of the linear group equals $GL(n,R)\suba\cap GL(n,R)\subb$.
Therefore, the natural map between both fundamental groups of $GL(n,R)$ is an isomorphism.

\

The homotopy of $A(G,\H)$ is locally trivial, in the sense that the local groupoids of $A$ are simply connected groupoids. This leads us to 
a more general family of atlases in which the map relating both groups $\pi_1^w(A,x)$ and $\pi_1(A,x)$ is an isomorphism.

\

\label{infimum}
Recall the definition of infimum from \cite{Bak1}, \cite{Bak2}. A groupoid atlas $A$ is {\it infimum} if for every local frame $s$ the set
$$\phi_s=\set{\a\in\phi_A\ |\ s\subset X\suba}$$
has an initial element. Note that an infimum groupoid atlas is, in particular, a good atlas.
This condition, when all local groupoids are simply connected, implies that both groups $\pi_1(A,x)$ and $\pi_1^w(A,x)$ are isomorphic by $p$.

\begin{thm}\label{pi1a equals pi1wa}
If $A$ is an infimum groupoid atlas such that every local groupoid $\G\suba$ of $A$ is simply connected, then the canonical map $p:\pi_1(A,x)\to\pi_1^w(A,x)$ is an isomorphism.
\end{thm}

\begin{proof}
Given a two elements local frame $\set{y,z}\subset X_A$, since $\phi_{\set{x,y}}$ has an initial element, 
there is a unique global arrow with source $y$ and target $z$. Then, $p$ maps bijectively the paths of global arrows into the paths of edges. 
The 2-simplices of $VA$ arise from the three elements local frames. Under this bijection, any three edges which are the faces of a 2-simplex in 
$VA$ correspond to a homotopy triangle.
\end{proof}

\begin{remark}
We propose a stronger version of the definition of infimum atlases (the definition of above is the original introduced in \cite{Bak1}) that 
seems to fit better from the strong point of view. We might say that $A$ is {\it infimum} if $\phi_s$ has an initial element for every simplex $s$ of $NA$.
When $A$ is such that all its local groupoids are simply connected, the original and the strong definition of infimum agree, since, in this case, a simplex is essentially determined by its underlying local frame.
\end{remark}

\begin{remark}
Theorem \ref{pi1a equals pi1wa} remains true when all the local groupoids are simply connected and $A$ satisfies the following condition: 
the sets $\phi_s$ are filtered for all local frame $s$.
Note that this condition is weaker than the infimum condition, but it is sufficient to prove the result.
\end{remark}

\

%Under these hypotheses, one can show that the weak and strong nerves of $A$ defined in the next section, are equal (\ref{na equals nwa}). The last proposition is, in fact, a consequence of that result.

%\

We finish this section considering the atlas $A=A(X,\U)$, where $X$ is a topological space and $\U$ is an open cover of $X$. The atlas $A$ is infimum if $\U$ is closed under intersections, and we have seen in \ref{topological example 2} that $A$ is irreducible if $U$ is path connected for all $U\in\U$.
When these conditions hold, the strong fundamental group of $A$ equals $\pi_1(X)$ (cf. \ref{fundamental group}) and the weak one is, by Dowker's theorem, the fundamental group of the nerve of the cover.

By last proposition, these groups are isomorphic when each $U$ is simply connected. By the remark of above, it is sufficient for $\U$ to be closed under finite intersections. Thus, we obtain as a corollary of \ref{pi1a equals pi1wa} an alternative proof of the following result.

\begin{corollary}
Let $X$ be a topological space and $\U$ an open cover by simply connected open subsets. If $\U$ is closed under finite intersections, then the group $\pi_1(X,x)$ equals the fundamental group of the nerve of the cover $N\U$.
\end{corollary}

\

\section{Nerves for Atlases}\label{s5}

Let $\ss$ the category of simplicial sets and simplicial morphisms.
We introduce two functorial constructions $\ga\to\ss$ for the \it nerve \rm of a groupoid atlas $A$.
The first one, $N^wA$, is based on the cover by components of local groupoids. The second one, $NA$, preserves more 
information about the local groupoids. Using these constructions, we can define homology theory of groupoid atlases and also a weak and a strong version of the classifying space of a groupoid atlas.

% classifying spaces $B^wA=|N^wA|$ and $BA=|NA|$.

\subsection{The weak nerve $N^w(A)$}

\begin{definition}
The {\it weak nerve} $N^w A$ of a groupoid atlas is the simplicial set whose $n$-simplices are the sequences of $n+1$ elements of the same component of a local groupoid,
$$(N^wA)_n=\set{(x_0,\dots,x_n)\ |\ \set{x_0,\dots,x_n}\text{ is an $\a$-frame for some $\a$}}$$
equipped with usual faces and degeneracies: the face $d_i$ erases the $i$-th element and the  degeneracy $s_j$ repeats the $j$-th.
\end{definition}

Given a weak map of groupoid atlases $f:A\to B$, we have $f_*:N^wA\to N^wB$ defined by $(x_0,\dots,x_k)\mapsto(f(x_0),\dots,f(x_k))$. 
It is well defined since $f$ preserves local frames. Note that this construction is functorial: $id_*=id$ and $(f\circ g)_*=f_*\circ g_*$.

\

Consider $a(\Delta[n])$ the $n$-simplex in the category $\ga$. Since it has only one index, the entire set $X_{a(\Delta[n])}=\set{0,...,n}$ is a local 
frame. A weak map $a(\Delta[n])\to A$ is a function $\set{0,\dots,n}\to X_A$ such that its image is a local frame. Thus, a simplex $s\in (N^wA)_n$ can be viewed as a weak map $a(\Delta[n])\to A$. Faces and degeneracies are, under this correspondence, compositions with the canonical inclusions $a(\Delta[n-1])\to a(\Delta[n])$ and projections $a(\Delta[n+1])\to a(\Delta[n])$, respectively. Hence,
%The nerve defined as before is the {\it singular functor} represented by $\underline{n}\mapsto a(\Delta[n])$, in the sense of \cite{FL}.
$$(N^w A)_n=Hom_{weak}(a(\Delta[n]),A).$$

The weak nerve of $A$ is the simplicial set that naturally arises from the complex $VA$ (cf. \cite{Sp}). The geometric realization $|N^wA|$ has the homotopy type of the polyhedron induced by $VA$, so it has the same fundamental group and the same homology groups.

\

\subsection{The (strong) nerve $N(A)$}

The weak nerve is constructed from the covering of $X_A$ by the components of each $\G\suba$. It 
has no more information about the groupoid structure than that. The strong version of the nerve appears naturally when one looks for a construction that preserves the local information.

\begin{definition}
Let $A$ be a regular atlas. The {\it nerve} $NA$ is the simplicial set whose $n$-simplices are
$$NA_n=\set{x_0\xto{g_1}x_1\xto{g_2}...\xto{g_n}x_n\ |\  g_1,...,g_n\in\G\suba\text{ for some }\a}/\sim$$
where we identify $(g_1,...,g_k)$ with $(\phiab(g_1),...,\phiab(g_n))$, its image through the structure functors $\phiab$. 
The face and degeneracy maps are defined as it is usual for nerves of categories: a face composes two arrows, a degeneracy inserts an 
identity (cf. \cite{Se}). They pass to the quotient since $\phiab$ are functors.
\end{definition}

\begin{remark}
Given  a groupoid $\G$, the nerve $N(a\G)$ coincides with the nerve $N\G\suba$ defined as usual for categories.
\end{remark}

Note that the $0$-simplices of an atlas $A$ are the points and the $1$-simplices are the global arrows. In general, the set $NA_k$ is a quotient of the disjoint union of the $k$-simplices of the nerves of the local groupoids $(N\G\suba)_k$.

\

Given $f:A\to B$ a map of groupoid atlases, it determines a map between the nerves $f_*:NA\to NB$ by the formula
$$f_*[(g_1,...,g_k)]=[(\G_f(\a)(g_1),...,\G_f(\a)(g_k))]$$
for  arrows $g_1,...,g_k$ in $\G\suba$. This definition does not depend on the representative of the class $[(g_1,...,g_k)]\in NA_k$.

\begin{proposition}\label{iso}
If $f\lhd f'$ then $f_*=f'_*$. Therefore, equivalent atlases have isomorphic nerves.
\end{proposition}

\begin{proof}
The simplex $[(g_1,...,g_k)]$ is mapped by $f_*$ and $f'_*$ into $[(\G_f(\a)(g_1),...,\G_f(\a)(g_k))]$ and $[(\G_{f'}(\a)(g_1),...,\G_{f'}(\a)(g_k))]$, respectively. Since $f\lhd f'$, we have $f(\a)\ls f'(\a)$ and $\phi_{\phi_f(\a)}^{\phi_{f'}(\a)}\circ\G_f(\a) = \G_{f'}(\a)$.
$$\xymatrix{
 & \G\suba \ar[dl]_{\G_f(\a)} \ar[dr]^{\G_{f'}(\a)} \\
\G_{\phi_f(\a)} \ar[rr]_{\phi_{\phi_f(\a)}^{\phi_{f'}(\a)}} & & \G_{\phi_{f'}(\a)}}$$
We conclude that $(\G_f(\a)(g_1),...,\G_f(\a)(g_k))\sim(\G_{f'}(\a)(g_1),...,\G_{f'}(\a)(g_k))$.
\end{proof}

The definition of the nerve $NA$ can be extended to non regular atlases, defining $NA$ as the nerve of its regularization:
$$NA=N(rA)$$

By the last proposition, this identity  remains true even when $A$ is already regular, since in this case,  $A\sim rA$ (cf. \ref{a sim ra}).

\

Recall that the $k$-simplices of the nerve $NC$ of a small category $C$ can be presented as $$NC_k=Hom_{\cat}(\set{0\to1\to...\to k},C),$$
and when $C$ is a groupoid, by the universal property of localization, we have
$$NC_k=Hom_{\gpd}(\Delta[k],C).$$

It is clear that $Hom_{\ga}(a(\Delta[k]),A)$ equals the disjoint union of the $k$-simplices of $N\G\suba$,
since a map $a(\Delta[k])\to A$ picks an index $\a$ of $\phi_A$ and gives a functor $\Delta[k]\to \G\suba$. Note that two maps $a(\Delta[k])\to A$ are equivalent 
if and only if the simplices that they define are identified in $NA$. Thus,
$$NA_k=Hom_{\ga}(a(\Delta[k]),A)/\sim\ =Hom_{[\ga]}(a(\Delta[k]),A).$$

Therefore, as it happens with the weak nerve, the functor $N:[\ga]\to \ss$ is an example of {\it singular functor} in the vocabulary of \cite{FL}. It means that there exists a functor $\theta:\Delta\to [\ga]$
from the category of finite ordinal numbers such that 
$NA_k=Hom_{[\ga]}(\theta(\set{0<...<k}),A)$,
and the faces and degeneracies are given by composition with $\theta(x)$, where $x$ is an elementary injection or surjection in $\Delta$. In this case, the functor $\theta$ is $\set{0<...<k}\mapsto a(\Delta[k])$. 
%The functor $\theta$ has an extension
%$$\Gamma_{\theta}:\ss\to [\ga],$$
%some kind of {\it groupoid atlas realization}, which is left adjoint to the nerve $N:[\ga]\to \ss$ (cf. \cite{FL}, proposition 2.3. See also \cite{McL},
%chapter X). In fact, it is not necessary to prove that $[\ga]$ is cocomplete: recall that $a$ is left adjoint to In particular, the nerve functor preserves limits.

\

\begin{remark}\label{na is not kan}
The nerve of a groupoid atlas $NA$ is not in general  a Kan complex. It might happen that two 1-simplices cannot be extended to another $1$-simplex: 
take local groupoids $\G\suba$ and $\G\subb$ and let $x,y$ be objects of $\G\suba$ and $\G\subb$, respectively. Let $z$ be an object shared 
by these two groupoids. Suppose that $\set{x,y,z}$ is not a local frame. If $g:x\to z$ and $g':z\to y$ are arrows of $\G\suba$ and $\G\subb$ respectively, then there is not a 2-simplex $s$ in $NA$ such that $d_0(s)=g'$ and $d_2(s)=g$.
\end{remark}

\begin{example}
In the nerve of the groupoid atlas corresponding to the global action $\A(D_3,\H)$ discussed in the introduction, there is no
 $s\in NA_2$ satisfying $d_0(s)=1\to r$ and $d_2(s)=r\to s\cdot r$.
\end{example}

\

If the atlas $A$ is irreducible, the set $NA_k$ is the colimit over $\phi_A$ of the sets $(N\G\suba)_k$, since it is the largest 
quotient of the disjoint union that makes  the diagrams
$$
\xymatrix{(N\G\suba)_k\ar[rr]^{\phiab}\ar[rd] & &(N\G\subb)_k\ar[dl]\\ & (\coprod (N\G\suba)_k)/\sim}
$$
commutative. Since limits and colimits in $\ss$ can be computed coordinatewise, we have
$$NA=\colim{\phi_A}{N\G\suba}.$$
When $A$ is not irreducible, there is not a well defined diagram in the category of sets since the functions $N\G\suba\to N\G\subb$ are partially defined.

\begin{remark}\label{na equals nia}
The nerve $NA$ of an atlas $A$ is equal to the nerve $N(iA)$ of the irreducible atlas $iA$, since both atlases
 have the same simplices by the universal property of $iA$. Then, the nerve of an arbitrary atlas can be computed as the nerve of an irreducible atlas. In the rest of this section we will assume, without loss of generality, that $A$ is irreducible.
\end{remark}

\medskip

\begin{proposition}\label{pi1 de na}
The fundamental group of the atlas $A$ is equal to $\pi_1(NA,x)$.
\end{proposition}

\begin{proof}
We have seen that the fundamental group of $A$ is the set of paths of global arrows modulo the relations generated by $(g,h)\sim(hg)$.
 Note that $g,h$ and $hg$ are the three faces of the 2-simplex $(g,h)$, so the group $\pi_1(A,x)$ is the set of paths of 1-simplices 
that start and end in $x$ modulo the simplicial homotopies, i.e. $\pi_1(A,x)$ equals $\pi_1(NA,x)$.
\end{proof}

\begin{corollary}
Let $A$ be a groupoid atlas and $x\in X_A$. The map $iA\to A$ induces an isomorphism $\pi_1(iA,x)\to\pi_1(A,x)$.
\end{corollary}

Last proposition also gives an alternative proof of \ref{pi1 y sim}.

\

\subsection{The colimit groupoid $\G(A)$}\label{colimit groupoid}

In  theorem \ref{fundamental group} we introduced another algebraic object related to $A$, namely 
the {\it colimit groupoid} over the diagram $\G_A$, which we denoted $\G(A)$.
$$\G(A)=\colim{\phi_A}{\G\suba}$$
Let $f=(\phi_f,\G_f):A\to B$ be a map in $\ga$. The family $\set{\G\suba\to \G_{\phi_f(\a)}\to \G(B)}_{\a\in\phi_A}$ commutes with 
the structural functors $\phi\suba^{\b}$, hence induces a functor $f_*:\G(A)\to \G(B)$ by the universal property of $\G(A)$.
With this definition on arrows, $\G:\ga\to\gpd$ becomes a left adjoint for the functor $a:\gpd\to\ga$.
$$Hom_{\gpd}(\G(A),G)\equiv Hom_{\ga}(A,aG)$$
The functor $\G$ factors through $[\ga]$, since a map and a corestriction of it give rise to the same family $\set{\G\suba\to \G_{\phi_f(\a)}\to \G(B)}_{\a\in\phi_A}$.

\

It is easy to see that $Obj\ \G(A)=X_A$. Recall that an arrow of $\G(A)$ is a path of colimit arrows modulo the smallest
equivalence class that contains the local identities and compositions. Here, a {\it colimit arrow} $x\to y$ means an element of $colim_{\phi_A}\ Hom_{\G\suba}(x,y)$. In other words, it is an arrow of the colimit graph of the underlying graphs of $\G\suba$.

\begin{remark}
Note that a global arrow of $A$ is the same as a colimit arrow. They are the arrows of $\G(A)$ that are in the image of some $\G\suba\to\G(A)$.
\end{remark}

The functors $i\suba:\G\suba\to\G(A)$ give rise to simplicial maps $Ni\suba:N\G\suba\to N\G(A)$. Last remark can be generalized to the following result.

\begin{proposition}
The nerve $NA$ is the subsimplicial set of $N\G(A)$ that consists of those simplices that are in the image of $Ni\suba$, for some $\a$.
\end{proposition}

\begin{proof}
Let $S$ be the union $\bigcup Ni\suba(N\G\suba)$. With the faces and degeneracies of $N\G(A)$, $S$ results a simplicial set. Recall that an $n$-simplex $s$ of $NA$ is the class of a functor
$$\Delta[n]\xto{s} \G\suba$$
for some $\a$, under the quotient map that identifies $s$ with $\phiab\circ s$ for all $\a\ls\b$.
We assign to $s$ the composition $i\suba\circ s:\Delta[n]\to \G(A)$. This factors through the quotient, so there is
 a function $NA_n\to S_n$. Since this function preserves faces and degeneracies, it gives rise to a simplicial map $NA\to S$. 
It is surjective by definition of $S$. It is not difficult to prove that it is also injective.
%$$\xymatrix{ & & \G\subb\\ \Delta[n]\ar[rru]^{s'}\ar[rrd]_s \\ & & \G\suba \ar[uu]^{\phiab}}$$
% We ca
%, and the $n$-simplices of $S$ can be identified with the maps
%$$\Delta[n]$$
\end{proof}

\begin{remark}\label{na neq nga}
It is well known that the nerve of a groupoid is a Kan complex, therefore the nerve $NA$ is not equal  in general to $N\G(A)$ (see \ref{na is not kan}).
\end{remark}

\

\subsection{The classifying space of a groupoid atlas}

Given an atlas $A$, we can associate to it a topological space $BA$ via its nerve. In this subsection we introduce this construction and
 relate it with  other spaces associated to $A$. 
 
\begin{definition}
The {\it classifying space} $BA$ of an atlas $A$ is the geometric realization of its nerve $NA$.
\end{definition}

Like the geometric realization of any simplicial set, $BA$ is a CW-complex with a cell for each non-degenerated simplex of $NA$. 
As a consequence of our results on the nerve of a groupoid atlas, we obtain the following propositions.
\begin{proposition}
Let $A$ be a groupoid atlas. Then the canonical map $\varphi_A:iA\to A$ induces a homeomorphism $B(\varphi_A):BiA\to BA$.
\end{proposition}
\begin{proposition}
Let $f,g:A\to B$ be two equivalent maps. The induced continuous functions $Bf$ and $Bg$ are equal. In particular, equivalent atlases have homeomorphic classifying spaces.
\end{proposition}

\begin{remark}
Given  a simplicial complex $K$, its associated groupoid atlas $a(K)$ is infimum and all its local groupoids are simply connected. 
Then, by proposition \ref{na equals nia} its weak nerve equals its strong nerve. It follows that $B(aK)=B^w(aK)$. Since $B^w(aK)$ is homotopy equivalent
 to the polyhedron $K$, we conclude that any {\it interesting homotopy type} can be obtained as the classifying space of a groupoid atlas.
\end{remark}

\begin{remark}
An algebraic loop in $A$ induces a topological one in $BA$, since $BL\cong\R$ and a map of groupoid atlases $L\to A$ which
 stabilizes can be restricted to a closed real interval. This assignation preserves homotopy classes and induces
 an isomorphism $\pi_1(A,x)\to\pi_1(BA,x)$, which is just the composition $\pi_1(A,x)\xto{\sim}\pi_1(NA,x)\xto{\sim}\pi_1(BA,x)$.
\end{remark}

\subsection*{Relation between the weak and the strong nerves}

\label{p en nervios}
Given  a groupoid atlas $A$, there is a canonical projection $p:NA \to N^wA$ given by
$$(x_0\xto{g_1}x_1\to...\xto{g_n}x_n)\mapsto (x_0,...,x_n)$$
Every simplex of $N^wA$ is in the image of $p$, thus $p:(N^wA)_k\to NA_k$ is onto for all $k$.

\begin{proposition}\label{na equals nwa}
If $A$ is infimum and all its local groupoids are simply connected, then $p:NA\to N^wA$ is an isomorphism.
\end{proposition}

\begin{proof}
The argument is analogous to the one used in $\ref{pi1a equals pi1wa}$. It is sufficient to prove that over each simplex of $N^wA$ there is one and only one simplex of $NA$.
\end{proof}

In general, this is not longer true. As we can see in the following example, $p$ might have no inverse, even no section.

\begin{example}\label{rombitos}
Let $X_A=\set{a,b,c,d}$, $\phi_A=\set{1,2}$ discrete. Let $\G_1$ and $\G_2$ be the simply connected groupoids with objects $\set{a,b,c}$ and $\set{a,b,d}$ respectively.

\

\centerline{\scalebox{1}{\includegraphics{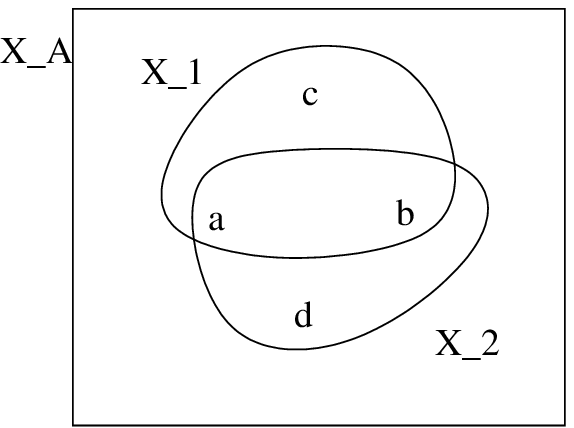}}}

\

Consider $(a,b)\in (N^wA)_1$. Suppose that there exists a section $i:N^wA\to NA$  for $p$. Then $p(i(a,b))=(a,b)$, and $i(a,b)=a\xto{1}b\text{ or }a\xto{2}b$, arrows of $\G_1$ and $\G_2$, respectively. Since $i(a,b,c)=a\xto{1}b\xto{1}c$, $i(a,b,d)=a\xto{2}b\xto{2}d$ and since $i$ commutes
with the face maps then we have
$$a\xto{1}b=d_2\circ i(a,b,c)=i\circ d_2(a,b,c)=i(a,b)= i\circ d_2(a,b,d)=d_2\circ i(a,b,d)=a\xto{2}b$$
which is a contradiction.

\end{example}

Compare the spaces $BA$ and $B^wA$.

\

\centerline{\scalebox{1}{\includegraphics{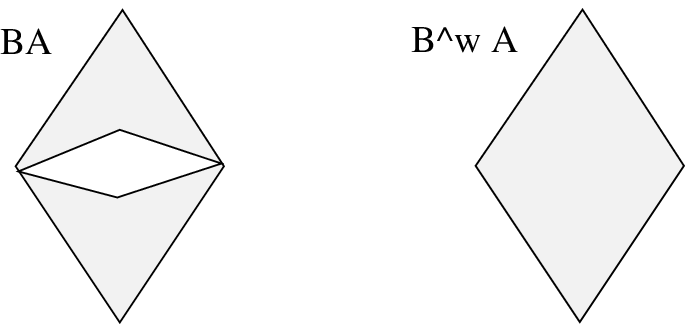}}}

\

In $BA$ one can notice the existence of a non trivial algebraic loop $a\xto{1}b\xto{2}a$, while the space $B^wA$ is simply connected.

\

\subsection{A little of homology}

\begin{definition}\label{def de homologia}
Let $A$ be a groupoid atlas and let $R$ be a commutative ring. The {\it homology of $A$ with coefficients in $R$} is the homology of the associated simplicial set $NA$.
$$H_n(A,R)=H_n(NA,R)$$
Explicitly, for each $n$ we put $C_n(A,R)=R[NA_n]$, the free $R$-module with basis $NA_n$, and the boundary map $d:C_n\to C_{n-1}$ is 
defined, as usual, by 
$$d(x)=\sum_i(-1)^id_i(x)$$
in the basis elements. We will simply denote $C_n(A)$ and $H_n(A)$ when $R=\Z$.
\end{definition}

%It seems interesting to try to relate this with topological spaces, or manifolds. Given an open cover, we can associate to it groupoid atlases equivalents to $\A(X,\U)$, with the same data but simplier than it, and compute the fundamental group, or even take some information about the homology groups from them.

\begin{generalities}

\

\begin{enumerate}

\item Two equivalent atlases have the same homology groups, since their nerves are isomorphics (cf. \ref{iso}).

\item Similarly, the homologies of $A$ and $iA$ are the same. Thus, we can restrict our attention to the homology of irreducible atlases.

\item If $\phi_A$ has a single element $\a$, the homology groups of $A$ are the homology groups of the groupoid $\G\suba$, which are the direct sum of the homologies of the vertex groups of each component of $\G\suba$.

\item If $n=0$, the group $H_0(A,R)$ is free with one generator for each component of $A$. Thus, we have an isomorphism
$$H_0(A,R)=R[\pi_0(A)].$$

\item When $R=\Z$, the groups $\pi_1(A)$ and $H_1(A)$ are related by a quotient map $\pi_1(A) \epi H_1(A)$ with kernel the commutator of the fundamental group (the first one coincides with $\pi_1(NA)$ and the second one with $H_1(NA)$). Thus, the classical equation
$$H_1(A)=\pi_1(A,a)/[\pi_1(A,a),\pi_1(A,a)]$$
remains true in this context.
\end{enumerate}
\end{generalities}

\begin{example}
Recall that a simplicial set and its geometrical realization share all their homology groups. We can 
compute the homology of the atlas $A$ of \ref{rombitos} from its classifying space, since it has the homotopy type of $S^1$.
$$H_n(A)=
\begin{cases}
 0 & \text{ if } n\neq 0,1\\
 \Z & \text{ if } n=0,1
\end{cases}$$
\end{example}

Given $A$ a groupoid atlas, note that the chain complex $C_*(A,R)$ is the colimit of the complexes $C_*(\G\suba,R)$: the free functor $R[-]$ is left adjoint to the forgetful functor, hence it preserves colimits.
$$C_*(A,R)=R[NA]=\colim{}{R[N\G\suba]}=\colim{}{C_*(\G\suba,R)}$$

\begin{example}\label{a,h y bH}
In the case of $A(G,\H)$, the chain complex $C_*(A)$ coincides with the chain complex $\b(\H)$, defined in \cite{AH} as the colimit of the $\Z[G]$-complexes induced by the non normalized homogeneous bar resolutions $\b(H\suba)$ of the groups $H\suba$. The homology groups of $A$ are the homology groups of the nerve of the cover of $G$ by the $H$-orbits, $H\in\H$.
\end{example}

\begin{remark}
The nerve $NA$ is a simplicial subset of $N\G(A)$. They may differ (\ref{na neq nga}) as the functor $N:\gpd\to\ss$ does not 
preserve colimits. When $A$ is a connected atlas, $H_*(\G(A))$ equals the homology groups of its fundamental group. We deduce that 
$H_*(A)\ncong H_*(\G(A))$ in general, since $BA$ can reach the homotopy type of any CW-complex, and the fundamental group of a CW-complex does not determine all its homology groups.
\end{remark}

\

\subsection{Relation between $HA$ and the local homology groups $H\G\suba$}

Consider the map $p=\coprod i\suba:\coprod_{\a\in\phi} N\G\suba\to NA$. It induces an epimorphism of chain complexes $p_*$, with kernel $Ker_*$. The short exact sequence 
$$0\to Ker_*\to \bigoplus_{\a\in\phi} C_*(\G\suba)\xto{p_*} C_*(A)\to 0$$
induces a long exact sequence relating the homologies of the complex $Ker_*$, the local groupoids $\G\suba$ and the groupoid atlas $A$.
$$...\xto{\partial} H_n(Ker)\to \bigoplus H_n(\G\suba)\to  H_n(A)\xto{\partial}  H_{n-1}(Ker)\to ...$$
In some cases, the complex $Ker_*$ can be written as a direct sum of some of the complexes $C_*(\G\suba)$. In these cases, the homology 
of $A$ can be computed from the homology of the local groupoids.

\begin{example}
Let $A$ be the {\it 1-sphere}, defined in \ref{sc example}. The set $X_A$ is $\set{0,1,2}$, the set of indices $\phi_A$ consists of 
 the proper  non empty subsets of $X_A$, ordered by  inclusion, and  for every $s$, the local groupoid $\G_s$ is the tree over $s$. 
Since it is infimum and all its local groupoids are simply connected, by (\ref{na equals nwa}) we have
$$NA_n=N^wA_n=\set{0,1}^{n+1}\cup\set{0,2}^{n+1}\cup\set{1,2}^{n+1}$$
For each $s$, the nerve of the local groupoid $\G_s$ is $(N\G_s)_n=s^{n+1}$. Then, there is a short exact sequence of complexes
$$0\to(C_*(\G_{\set0})\oplus C_*(\G_{\set1})\oplus C_*(\G_{\set2}))^2\xto{i_*}\bigoplus_{s\in\phi_A}C_*(\G_s)\xto{p_*} C_*(A)\to 0$$
where
$$i_*(\sum c_i\set i^{n+1},\sum c_i'\set i^{n+1})_s= \begin{cases}
c_i+c_i' &\ if\ s=\set i\\
-c_i-c'_j &\ if\ s=\set{i,j},\ j\equiv i+1 \text{ mod }3
\end{cases}$$
Since the groupoids $\G_s$ are simply connected, the non trivial part of the long exact sequence relating the homology groups is, in this case,
$$0\to H_1(A)\to \Z^6\to \Z^6 \to H_0(A)\to 0.$$
The group $H_0(A)$ is isomorphic to $\Z$ because $A$ is connected, and so $H_1(A)\cong \Z$: it is free of rank 1 since $H_1(A)\mono\Z^6$ and  the 
Euler characteristic of the last complex is 0.

% The projection $q:\set{0,1}^{n+1}\sqcup\set{0,2}^{n+1}\sqcup\set{1,2}^{n+1}\to NA_n$ is the coequalizer of the canonical inclusions $i_0,i_1:\set{0}^{n+1}\sqcup\set{1}^{n+1}\sqcup\set{2}^{n+1}\rightrightarrows\set{0,1}^{n+1}\sqcup\set{0,2}^{n+1}\sqcup\set{1,2}^{n+1}$.  Thus, applying the free functor to this arrows we obtain short exact sequences
% $$0\to \bigoplus_{s=\set0,\set1,\set2}\Z[s^{n+1}]\xto{(i_0)_*-(i_1)_*}\bigoplus_{s=\set{0,1},\set{1,2},\set{0,2}}\Z[s^{n+1}]\xto{q_*} C_n(A)\to 0$$
% The nerve of the local groupoid $\G_s$ is $(N\G_s)_n=s^{n+1}$ for $\emptyset\neq s\subsetneq X_A$. Then, there is a short exact sequence of complexes
% $$0\to\bigoplus_{s=\set0,\set1,\set2}C_*(\G_s)\to\bigoplus_{s=\set{0,1},\set{1,2},\set{0,2}}C_*(\G_s)\to C_n(A)\to0$$
% The map $q_*$ is a restriction of $p_*$. In this case, $Ker(p_*)=Ker(q_*)\oplus Ker(q_*)$.
\end{example}

\begin{remark}
Of course, the homology of $a(\partial\Delta[1])$ can be easily computed noting that its classifying space is homeomorphic to $S^1$, but we 
proceeded in this way to illustrate the general idea (see below).
\end{remark}

\begin{remark}
All the local groupoids are, in this example, simply connected. Then, by the long exact sequence, the homology groups $H_n(A)$ should be trivial 
for $n\gs 2$. This does not happen, for example, with the $n$-sphere ($n\gs 2$), whose local groupoids are also simply connected. 
Therefore, for many examples we cannot find any expression of $Ker_*$ as a sum of complexes $C_*(\G\suba)$.
\end{remark}
\smallskip

In the example of above, we use the usual presentation of the colimit $C_*(A)$ as the cokernel of the map
$$\bigoplus_{\a<\b\in\phi} C_*(\G\suba)\xto{j}\bigoplus_{\a\in\phi} C_*(\G\suba)$$
which is defined on the basis elements by $j(s,\a<\b)=(s,\a)-(\phiab(s),\b)$. 
Here  $(s,i)$ denotes the element of the factor $i$ given by the simplex $s$. In general $j$ is not a monomorphism, then 
one cannot identify $Ker_*$ with the sum $\bigoplus_{\a<\b\in\phi} C_*(\G\suba)$.

\begin{proposition}
The map $j$ is mono if and only if $\phi_s$ has no cycles (viewed as a graph) for all simplices $s$.
\end{proposition}
\begin{proof}
Note that $Ker(j)=\bigoplus Ker(j_n)$, where $j_n$ is the $n$-th component of $j$.

Let $k:\bigoplus_{\a\in\phi}C_n(\G\suba)\to\bigoplus_{\a\in\phi}C_n(A)$ be the map defined in the basis elements by 
$(s,\a)\mapsto (\-s,\a)$. Here $\-s$ denotes the class of the simplex $s$ in $NA$. Let $\sum_i c_i(s_i,\a_i<\b_i)$ be an element of $Ker(j_n)$.
We have
\begin{align*}
0&=k(j(\sum_i c_i(s_i,\a_i<\b_i)))\\
&=k(\sum_i c_i j(s_i,\a_i<\b_i))\\
&=k(\sum_i c_i(s_i,\a_i)-c_i(\phiab(s_i),\b_i))\\
&=\sum_i c_ik(s_i,\a_i)-c_ik(\phiab(s_i),\b_i)\\
&=\sum_i c_i(\-{s_i},\a_i)-c_i(\-{s_i},\b_i)\\
\end{align*}
Then, for all $s\in NA_n$ we have $0=\sum_{\-{s_i}=s} c_i(s,\a_i)-c_i(s,\b_i)$ and $d(\sum_{\-{s_i}=s} c_i(\a_i,\b_i))=0$, where $d$ is the boundary map of the simplicial chain complex of $\phi_s$. 
This way, a non trivial element of $Ker(j_n)$ gives a cycle in some $\phi_s$. The converse is clear.
\end{proof}

\begin{remarks}\

\begin{itemize}
\item If there exist $\a,\b,\c\in\phi_A$ such that $\a<\b<\c$, then $j$ is not a monomorphism.
\item $j$ could be a monomorphism even if $\phi_A$ has non trivial cycles, e.g. $A$ the 1-sphere.
\end{itemize}
\end{remarks}

We expose a simple aplication of last result. Let $A$ be a groupoid atlas such that $\phi_A$ is discrete ($\a\ls\b\then\a=\b$). Recall that, since $A$ is not good, $NA=NrA$ by definition, 
and note that $(\phi_{rA})(s)$ is discrete for $s$ of dimension $\geq 1$ and is a star for $s$ a point. Thus, we have

\begin{corollary}
If $\phi_A$ is discrete, then $H_n(A)=\bigoplus\suba H_n(\G\suba)$ for $n\gs 2$. For $n=1$ we have $H_1(A)=(\bigoplus\suba H_1(\G\suba))\oplus F$ with $F$ a free abelian group.
\end{corollary}

\smallskip

The procedure of above can be emulated in other contexts, considering a cofinal subset of local groupoids $\set{\G\suba}_{\a\in S}$ to make 
the canonical map $p:\bigoplus_{\a\in S}C_*(\G\suba)\to C_*(A)$ an epimorphism (when $\phi_A$ is finite, then $S$ could be the set of maximal objects), 
and then
construct $Ker_*(p)$ as the sum of the chain complexes associated to certain local groupoids.

\email{ mdelhoyo$@$dm.uba.ar, gminian$@$dm.uba.ar}

\begin{thebibliography}{xxxx}

\bibitem{AH}
H. Abels, S. Holz.
\newblock \em Higher generation by subgroups.
\newblock \em J.Algebra 160, 311--341, 1993.

\bibitem{Bak1}
A. Bak.
\newblock \em Global Actions: The algebraic counterpart of a topological space.
\newblock \em Russian Math Surveys, 5; 955--996, 1997.

\bibitem{Bak2}
A. Bak.
\newblock \em Topological methods in algebra.
\newblock \em Rings, Hopf Algebras and Brauer Groups, Lec. Notes in Pure and Appl. Math., 197, 1998.

\bibitem{BBMP}
A. Bak, R. Brown, E.G. Minian, T. Porter.
\newblock \em Global actions, groupoid atlases and applications.
\newblock \em Journal of Homotopy and Related Structures 1(2006) 101-167.



\bibitem{Bass}
H.Bass.
\newblock \em Algebraic K-Theory.
\newblock \em Benjamin, 1968.

\bibitem{Brown}
R. Brown.
\newblock \em Topology and Groupoids.
\newblock \em Booksurge, 2006. 

\bibitem{FL}
R. Fritsch, D. M. Latch.
\newblock \em Homotopy Inverses for Nerve.
\newblock \em Mathematische Zeitschrift, 1981.

\bibitem{GZ}
P. Gabriel, M. Zisman.
\newblock \em Calculus of Fractions and Homotopy Theory.
\newblock \em Springer-Verlag, 1967.

\bibitem{McL}
S. Mac Lane.
\newblock \em Categories for the Working Mathematician.
\newblock \em Springer, 1971.

\bibitem{May}
J.P. May.
\newblock \em A Concise Course in Algebraic Topology.
\newblock \em Chicago Lectures in Mathematics, 1999.

\bibitem{Mil1}
J. Milnor.
\newblock \em Introduction to algebraic K-Theory.
\newblock \em Annals of Mathematics Studies, 72, Princeton, 1971.

\bibitem{Min1}
E.G. Minian. 
\newblock  \em Generalized cofibration categories and global actions.
\newblock  \em K-Theory  20, 37-95, 2000.

\bibitem{Min2}
E.G. Minian.
\newblock \em Lambda-cofibration categories and the homotopy categories of 
Global Actions and Simplicial Complexes.
\newblock \em Applied Categorical Structures 10; 1--21, 2002.

%\bibitem{Min4}
%E.G. Minian.
%\newblock \em Atlases of topological categories, local pospaces and dihomotopy.
%\newblock \em Preprint, 2005.

\bibitem{Por1}
T. Porter.
\newblock \em Geometric aspects of multiagent systems.
\newblock \em Electron. Notes Theoret. Comput. Sci. 81 (2003)


\bibitem{Qu}
D. Quillen.
\newblock \em Higher algebraic K-theory I.
\newblock \em Lectures Notes in Math 341;  85-147,\rm Springer 1973.

\bibitem{Se}
G. Segal.
\newblock \em Classifying spaces and spectral sequence.
\newblock \em Inst. Hautes \'Etudes Sci. Publ. Math. 34; 105-112, \rm 1968.


\bibitem{Sp}
E. Spanier.
\newblock \em Algebraic Topology.
\newblock \em Springer, 1966.


\bibitem{Swan}
R.G. Swan.
\newblock\em Algebraic K-Theory.
\newblock\em Lecture Notes in Mathematics, 76; Springer, 1968.

\end{thebibliography}
\end{document}